\documentclass[11pt]{article}
\usepackage{amsmath}
\usepackage{amssymb}
\usepackage{amsthm}
\usepackage{graphicx}
\newtheorem{Theorem}{Theorem}[section]

\newtheorem{example}[Theorem]{Example}

\newtheorem{remark}[Theorem]{Remark}

\newenvironment{Remark}{\begin{remark}\rm}{\end{remark}}
\newtheorem{Definition}[Theorem]{Definition}
\newtheorem{Algorithm}[Theorem]{Algorithm}

\newenvironment{Proof}
{\bigskip\noindent{\bf Proof.}~\rm}{\bigskip}
\begin{document}
\title{\huge Indefinite Integration Operator Identities and their
  Polynomial Approximations}

\author{{Frank Stenger}
\\{\sl Manager, SINC, LLC}
\\{\sl School of Computing}
\\{\sl Department of Mathematics}
\\{\sl University of Utah}
\\{\sl Salt Lake City, UT \ \ 84112 USA}}



\maketitle

\begin{abstract}
  \noindent The integration operators (*) $({\mathcal J}^+\,g)(x) =
  \int_a^x g(t) \, dt$\,, $({\mathcal J}^-\,g)(x) = \int_x^b g(t) \,
  dt$ defined on an interval $(a,b) \subseteq {\mathbf R}$ yield new
  identities for indefinite convolutions, control theory, Laplace and
  Fourier transform inversion, solution of differential equations, and
  solution of the classical Wiener--Hopf integral equations.  These
  identities are are expressed in terms of ${\mathcal J}^\mp$\, and
  they are thus esoteric.  However the integrals (*) can be
  approximated in many ways, yielding novel and very accurate methods
  of approximating all of the above listed relations.  Several
  examples are presented, using Legendre polynomial as approximations,
  and references are given for approximation of some of the operations
  using Sinc methods.  These examples illustrate for a class of
  sampled statistical models, the possibility of reconstructing models
  much more efficiently than by the usual slow Monte--Carlo
  $({\mathcal O}(N^{-1/2})$ rate.  Our examples illustrate that we
  need only sample at $5$ points to get a representation of a model
  that is uniformly accurate to nearly $3$ significant figure
  accuracy.
  \end{abstract}

\vspace{.2in} \noindent {{\bf Keywords}: Indefinite integration, indefinite
  convolution, Fourier transform inversion, Laplace transform
  inversion, Wiener-Hopf, differential equations,
  approximations}

\section{Introduction and Summary}
\setcounter{equation}{0}

\vspace{.1in} \noindent This paper presents some symbolic--like
approximations gotten from identities -- some previously known, and
some new -- of indefinite integration operators.  These operators are
defined for an interval $(a,b) \subseteq {\mathbf R}$ by the equations

\begin{equation}
  \begin{array}{rcl}
    ({\mathcal J}^+\,g)(x) & = & \displaystyle \int_a^x g(t) \, dt \\
    && \\
    ({\mathcal J}^-\,g)(x) & = & \displaystyle \int_x^b g(t) \, dt\,.
    \label{equation:1.1}
    \end{array}
\end{equation}
                                                                  
\noindent Here $g$ is a function defined on $(a,b)$\,, and we assume,
of course, that the integrals exist for all $x \in (a,b)$\,.  Combined
with convolution [Sp], Laplace and Fourier transform inversion, these
operators enable many new expressions for one--dimensional models, of:
control theory, for Laplace and Fourier transform inversion, for
solving ODE [DS] and PDE [Fc, Sc], and for solving Wiener--Hopf
integral equations.  These operations can be readily combined to get
multidimensional approximations.  And while these formulas, expressed
in terms of ${\mathcal J}^\mp$ are esoteric, and seemingly devoid of
any practical value, when replaced with certain types of approximation
in terms of interpolation on $(a,b)$\,, yield novel, accurate and
efficient methods of approximation.  The mode or operations of this
paper also enables defining models accurately in terms of statistical
samples, but instead of sampling over the whole interval, it suffices
to sample at just a small number of points on the interval\,.  Thus,
with respect to our illustrative one dimensional examples in this
paper, for which we sample at just 5 points, we can recover the model
using $5^2 = 25$ points, in $2$ dimensions, $5^3 = 125$ points in $3$
dimensions, etc.

\vspace{.2in} \noindent A formula for the approximation of the
indefinite integral was first published in [Si], albeit without proof.
The first proof was published by Kearfott [K] and later, two proofs
were published by Haber [Ha].  The author later presented a proof of
this result in \S4.5 of [Sp].  These approximation formulas were
applied in [Fc] and [Sc] to solve partial differential equations
(PDE), via use of explicit Laplace transforms derived by the author of
the Poisson, the heat and the wave equation Green's functions in one
and in more than one dimension.  In [STB], the Laplace transform of
the heat equation Green's function in ${\mathbf R}^3 \times
(0,\infty)$ was used to obtain a numerical solution of the
Navier--Stokes equations, and in [SKB] the Fourier transform of this
same Green's function was used to obtain a numerical solution of the
Schr\"odinger equation in ${\mathbf R}^3 \times (0,T)$.  More
recently, in [Srh], the author constructed a proof of the Riemann
hypothesis by use of the operators ${\mathcal J}^\mp$\,.

\vspace{.2in} \noindent In \S2 we present the operators ${\mathcal
  J}^\mp$\,, along with some of their properties, as well as methods
of approximating the operations of these operators.  In addition, we
present identities for optimal control, Laplace transform inversion,
and the solution of Wiener--Hopf integral equations, as well as
identities that are based on Fourier transforms for optimal control,
for Fourier transform inversion, for the solution of ordinary
differential equations, and for the solution of Wiener--Hopf integral
equations.  Some of these formulas -- the ones involving Laplace
transforms -- were previously known [Sp], whereas those related to
Fourier transforms are new.  We have also omitted the relation of
these operators with solution of PDE, since this aspect was covered
extensively in [Sc].

\vspace{.2in} \noindent In \S3 we illustrate explicit the application
of the esoteric formulas developed in \S2; the replacement of
${\mathcal J}^\mp$ of the formulas of \S2 with explicitly defined
matrices $A^\mp$ transform these esoteric identities into accurate and
efficient novel methods of approximation.  The matrices that have been
used to date are based on either Sinc or Fourier polynomials. Sinc
methods have previously been used to define these matrices (see e.g.,
[Sc], [STB], [SKB], [SA] and [SG]); we have restricted our examples to
using Legendre polynomials to define and use the matrices $A^\mp$\,,
since the use of other methods of approximation is similar.

\vspace{.2in} \noindent An important property of the matrices $A^\mp$
which enables the functions of the matrices $F(A^\mp)$ which are
gotten by replacement of ${\mathcal J}^\mp$ with $A^\mp$ in the
operator expression $F({\mathcal J}^\mp)$ to be well defined is that
the eigenvalues of $A^\mp$ are located on the right half of the
complex plane ${\mathbf C}$\,.  This was a 20--year conjecture for
Sinc methods; a proof of this conjecture was first achieved by Han \&
Xu [HX].  A proof was obtained by Gautschi \& Hairer [GH] for Legendre
polynomials, but a proof for the case of polynomials that are
orthogonal over an interval with respect to a positive weight function
is still an open problem; this author author offers \$300 for the
first proof or disproof, that the real parts of all of the eigenvalues
of the corresponding integration matrices defined in Definition 3.1 of
this paper have positive real parts.

\section{The Hilbert space and the operators}
\setcounter{equation}{0}

It is most convenient to work with operators in the setting of a
Hilbert space.  Let $(a,b) \subseteq {\mathbf R}$\,.  Our Hilbert
space is just the well--known space ${\mathbf H} = {\mathbf
  L}^2(0,b)$\,, with $b \in {\mathbf R}^+$\,.

\subsection{The operators ${\mathcal J}$}

Let the operators ${\mathcal J}^\mp$ be defined for ${\mathbf H}$\,.
The inverses of these operators have the property: if $G = {\mathcal
  J}^\mp\,g$\,, then $g = \left({\mathcal J}^\mp\right)^{-1}\,G$\,,
i.e., $g(t) = \mp \, \frac{d}{dt}\,G(t)$\,, whenever the derivatives
exist.

\subsection{Numerical ranges}

We mention here some properties of numerical ranges for the operators
${\mathcal J}^\mp$\,.

\begin{Definition}
  \label{definition:2.1}
  Let ${\mathbf H}$ be defined as above, and let the
  operators ${\mathcal J}^\mp$ be defined as in (1.1).    The 
  numerical range ${\mathbf W}$ of ${\mathcal J}^\mp$ in 
  ${\mathbf H}$  is defined as

  \begin{equation}
{\mathbf W} = \{({\mathcal J}^\mp\,f,f) : f \in {\mathbf H}\,, \ \ {\rm with}
\ \ \ \|f\| = (f,f)^{1/2} = 1\}\,.
\label{equation:2.1}
\end{equation}

\end{Definition}

\noindent ``Numerical range'' is synonymous with ``field of values''
with the latter being used more often for matrices.  The closure of
the numerical range of ${\mathcal J}^\mp$ contains the spectrum of
${\mathcal J}^\mp$\,.  Other properties of the numerical range can be
found, for example, in [GR] and in [Sh]\,.

\begin{Theorem}
     \label{theorem:2.2}
     Let ${\mathcal J}^\mp$ be defined as in Definition 2.1.  Then:

     (i.) The numerical ranges of ${\mathcal J}^\mp$ are contained in the
     closed right half plane, $\{z \in {\mathbf C} : \Re z \ge 0\}$\,.

     (ii) The real part of the numerical range of $\mp \left({\mathcal
       J}^\mp\right)^{-1}$ is $(1/2)(|f(b)|^2 - |f(a)|^2)$\,.
\end{Theorem}

\begin{Proof} Part(i.): We give a proof of the (i.)--part of this theorem 
  only for the case of ${\mathcal J}^+$\,, inasmuch as the proof for
  ${\mathcal J}^-$ is similar.

 \vspace{.1in} \noindent Let $g \in {\mathbf H}$ denote a
 complex valued function.  Then

\begin{equation}
  \begin{array}{rcl} 
    \Re \left({\mathcal J}^+\,g,g\right) 
    & = & \displaystyle \Re \int_a^b \overline{\int_a^x g(t)\,dt} \,
    g(x)\,dx\,,  \\
    && \\
   & = & \displaystyle \frac{1}{2} \, \left|\int_a^b g(x) \, dx\right|^2
   \geq 0\,, 
  \end{array}
  \label{equation:2.2}
\end{equation}

\noindent so that the inner product $({\mathcal J}^+\,g,g)$ is contained in
the closed right half plane.

\vspace{.1in} \noindent Part (ii.) We also omit the straight--forward
proof of this part of the theorem.
 \end{Proof} $\blacksquare$

\begin{Theorem}
\label{theorem:2.3}
  Let $(0,b) \subset {\mathbf R}$ be a finite interval,
  let the Hilbert space ${\mathbf H}$ be ${\mathbf
    L}^2(0,b)$\,.  Then  

  \begin{equation}
    \|{\mathcal J}^+\| \leq \frac{b}{\sqrt{2}}\,.
    \label{equation:2.3}
  \end{equation}

  \end{Theorem}

  \begin{Proof} Let $g \in {\mathbf L}^2(a,b)$\,.  Then we have

    \begin{equation}
\begin{array}{rcl}
\|{\mathcal J}^+\,g\|^2 & = & \displaystyle \int_0^b
  \overline{\int_0^x g(t)\,dt} \, \int_0^x g(\tau) \,
  d\tau \,  dx \\
 && \\
& \leq & \displaystyle \int_0^b \, \left(\int_a^x 1\,|g(t)|\,dt\right)^2 dx
  \\ 
&& \\
& \leq &  \displaystyle \int_0^b x \, dx \, \int_0^x|g(t)|^2\,dt \\ 
&& \\ 
& = & \displaystyle \frac{b^2}{2}\,\|g\|^2\,,
\end{array}
\label{equation:2.4}
    \end{equation}

\noindent which yields (2.3).
\end{Proof} $\blacksquare$

 \subsection{Indefinite Convolution via Fourier Transforms} 
 
 \vspace{.2in} \noindent We now take $(0,b) = {\mathbf R}^+$\,, where
 ${\mathbf R}^+$ denotes the interval $(0,\infty)$\,, we take $f^\mp
 \in {\mathbf H} = {\mathbf L}^2({\mathbf R}^+)$\,, and we assume the
 usual Fourier and inverse Fourier transforms defined by 

 \begin{equation}
   \begin{array}{rcl}
     \widehat{f^\mp}(y) & = & \displaystyle \int_{{\mathbf R}^+} f^\mp(x) \,
     e^{\mp\,i\,x\,y} \, dx \\
     && \\
 f^\mp(x) & = & \displaystyle \frac{1}{2\,\pi} \int_{\mathbf R}
 e^{\mp\,i\,x\,y} \, \widehat{f^\mp}(y) \, dy\,.
   \end{array}
   \label{equation:2.5}
 \end{equation}

 \newpage
\begin{Theorem}
  \label{Theorem:2.4}
  If $F^\mp := \widehat{f^\mp}$ is the Fourier transform of $f^\mp \in
  {\mathbf H}$\,, if ${\mathcal J}^\mp$ is supported on ${\mathbf
    R}^+$\,, if $(F^\mp)^\prime$ exists on $[0,\infty]$\,, and if
  $f^\mp$ is real--valued on ${\mathbf R}^+$\,, then
  
  \begin{equation}
    \begin{array}{rcl}
    (\mp i/{\mathcal J}^\mp\,F^\mp,F^\mp) & = & \pm \,  \displaystyle \pi \,
      \int_{{\mathbf R}^+}y\,|f^\mp(y)|^2 + \frac{i}{2} \,
      \left(\int_{{\mathbf R}^+} f^\mp(y) \, dy\right)^2\,, \\
    && \\
    (\mp i/{\mathcal J}^\mp\,F^\pm,F^\pm) & = & \mp \displaystyle \pi \,
      \int_{{\mathbf R}^+}y\,|f^\mp(y)|^2 - \frac{i}{2} \,
      \left(\int_{{\mathbf R}^+} f^\mp(y) \, dy\right)^2\,.
      \end{array}
    \label{equation:2.6}
  \end{equation}

  \end{Theorem}

\begin{Proof}  We shall prove this lemma only for the case of
  ${\mathcal J}^+$ acting on $F^+$\,, since the proofs for the other
  cases are similar.

\vspace{.1in} \noindent By taking $F^+ = u + i\,v$ where $F^+$ is the Fourier
  transform of $f^+ \in {\mathbf H}$ taken over ${\mathbf R}^+$\,, we
  have, 

\begin{equation}
 \begin{array}{rcl}
 \Im \left(i/{\mathcal J}^+\, F^+\,,F^+\right) & = &   
\displaystyle \Im \int_{{\mathbf R}^+} \overline{i\,(F^+)^\prime}(x) \,
F^+(x) \, dx \\ 
&& \\
& = &  \displaystyle  - \, \frac{1}{2} \, \left.(u^2(x) +
v^2(x))\right|_0^\infty\,, \\
&& \\
& = &  \displaystyle  \frac{1}{2} \, |F^+(0)|^2 \\
&& \\
& = &  \displaystyle \frac{1}{2} \, \left(\int_{{\mathbf R}^+} f^+(y) \,
    dy\right)^2\,. 
  \end{array}
  \label{equation:2.7}
\end{equation}

\noindent Note, the term in the last line cannot be negative.  If the
functions $f \in {\mathbf H}$ belongs to ${\mathbf L}^1({\mathbf
  R})$\,, then the term $u^2(\infty) + v^2(\infty) = 0$\,, since then
$|F(x)|^2 \to 0$ as $x \to \infty$\,, by the Riemann--Lebesgue lemma.

\vspace{.2in} \noindent Next, for the real part:

\begin{equation}
  \begin{array}{l} 
    \Re (i/{\mathcal J}^+\,F^+,F^+) = \\
    \\
    \displaystyle \Re \, -i \, \int_{{\mathbf R}^+}
    \overline{\frac{{\partial}}{{\partial} x} \, 
      \int_{{\mathbf R}^+} e^{i\,x\,y} \, f^+(y) \, dy} \,
    \int_{{\mathbf R}^+} e^{i\,x\,\eta} \, f^+(\eta) \, d\eta \, dx \\
  \\
  =  \displaystyle \Re -i \lim_{{\varepsilon} \to 0^+} \int_{{\mathbf
      R}^+} -i\,y\,\overline{f^+(y)} \, 
  \int_{{\mathbf R}^+} f^+(\eta) \, \int_{{\mathbf R}^+} 
  e^{-\,x({\varepsilon} + i\,(y-\eta))} \, dx \, d\eta \, dy \\
\\
=  \displaystyle \int_{{\mathbf R}^+} - \,  \int_{{\mathbf
      R}^+}  y\,\overline{f^+(y)} \, \lim_{{\varepsilon} \to 0} \int
  \frac{f^+(\eta)\,\varepsilon}{(x - \eta)^2 + {\varepsilon}^2} \, d\eta \, dy \\
\\
  = \displaystyle - \, \pi \, \int_{{\mathbf R}^+} y\,|f(y)|^2\,dy\,.
  \end{array}
  \label{equation:2.8}
\end{equation}

 \noindent The interchange of the order of integration in (2.8) is permitted
 since both functions $f(y)$ and $y\,f(y)$ belong to ${\bf
   L}^1({\mathbf R}^+)$\,.  

\vspace{.1in} \noindent This completes the proof of Theorem 2.4\,.
\end{Proof} $\blacksquare$

\subsection{Optimal Control}

Our model indefinite integrals for $x \in (0,b) \subseteq
(0,\infty)$ corresponding to given functions $f^\mp$ and $g$ take the
form

 \begin{equation}
 \begin{array}{rcl}
   q^+(x) & = & \displaystyle  \int_0^x f^+(x-t)\,g(t)\,dt, \\
   && \\
   q^-(x) & = & \displaystyle  \int_x^b f^-(x-t)\,g(t)\,dt.
 \end{array}
 \label{equation:2.9}
 \end{equation}

\noindent Given one or both functions $\widehat{f^\mp}$ in (2.5), we shall 
obtain a formula for determining $q^\mp$ on $(0,b)$\,, under the
assumption that the indefinite integration operators ${\mathcal
  J}^\mp$ are supported on $(0,b)$\,. 

\vspace{.1in} \noindent Novel explicit evaluations of the integrals
(2.9) were first obtained in [Sp], \S4.6, by use of Laplace
transforms.  Indeed, many new results were obtained using those
formulas, including novel explicit formulas for Laplace transform
inversions, novel explicit formulas for evaluating Hilbert transforms
(discovered independently by Yamamoto [Y] and Stenger [Sc], \S1.5.12),
and novel formulas for solving partial differential and convolution
type integral equations [Sc].  Included with each of these formulas
are very efficient and accurate methods of approximation -- most of
which are given in [Sc]; these usually are orders of magnitude more
efficient than the current popular methods of solving such equations.

\vspace{.1in} \noindent We shall now derive similar one dimensional
novel convolution formulas based on Fourier transforms.

  \begin{Theorem}
    \label{theorem:2.5}
Let ${\mathcal J}^\mp$ be defined as above, and have support on $(0,b)
\subseteq {\mathbf R}^+$\,, let the functions $f^\mp$ of equation
(2.6) belong to ${\mathbf H}$\,.  Let $q^\mp$ be defined as in
(2.6).  Then
 
 \begin{equation}
 q^\mp = \widehat{f^\mp}\left(\mp i\,/{\mathcal J}^\mp\right)\,g\,.
\label{equation:2.10}
 \end{equation}

 \end{Theorem}

 \begin{Proof} The proof resembles that given in \S 4.6 of [Sp] and
   \S1.5.9 of [Sc] for the case of Laplace transforms. We consider
   only the case of $f^+$\,, since the proof for the case of $f^-$ is
   similar.     The proof makes use of the following:
 \begin{itemize}
   \item $f^+$ has compact support $[0,b]$\,, -- inasmuch
       extension to $(0,\beta)$\,, or to an infinite interval can be
       carried out via the usual well--known procedure of analysis,
       and moreover, our assumptions are consistent with this possibility;
   \item By Theorem 2.4, the spectrum of $i/{\mathcal J}^+$ lies in
     the closed second quadrant of the complex plane, ${\mathbf C}$\,.
   \item As can be shown via an easy to prove inequality -- see
       [Sp], \S4.6 -- namely, that $\left\|{\mathcal
         J}^\mp\right\|_p \leq \alpha$ for all $p \in [1,\infty]$\,,
       where $(0\,,\alpha) \subset {\mathbf R}^+$ is the support of
       ${\mathcal J}^\mp$\,;
   \item By inspection of (2.6), $\widehat{f^+}$ is
     analytic in the upper half of the complex plane, and
      $\widehat{f^-}$ is analytic in the lower half; and
   \item The well--known formulas (2.6), as well as  

   \begin{equation} 
     \begin{array}{rcl}
          \displaystyle \left(\left({\mathcal J}^+\right)^n\,g\right)(x) & = & 
          \displaystyle \int_0^x 
          \frac{(x-t)^{n-1}}{(n-1)\,!} \, g(t) \, dt\,, \ \ \ n =
          1\,,2\,,\ldots\,,\\ 
           \\
          \widehat{f^\mp}(z)  & = & \displaystyle \mp
          \frac{1}{2\,\pi\,i} \int_{\mathbf R} 
          \frac{\widehat{f^\mp}(t)}{t - z} dt\,, \ \ \ \mp\,\Im z > 0\,.
     \end{array}
    \label{equation:2.11}
   \end{equation}
   
 \end{itemize}       

\noindent Hence, applying the above points, for $x \in (0,b)$\,, and with $I$
denoting the identity operator, we get

\begin{equation}
  \begin{array}{rcl}
q^+(x) & = & \displaystyle \int_0^x f^+(x-\xi) \, g(\xi) \, d\xi \\
    && \\
    & = & \displaystyle \int_0^x \,\left(\frac{1}{2\,\pi} \, \int_{\mathbf R}
\widehat{f^+}(t) \, e^{-\,i\,t\,(x - \xi)}\,dt \,
\right)\,g(\xi)\,d\xi \\  
&& \\
& = & \displaystyle \frac{1}{2\,\pi} \,\int_{\mathbf R} \left(\int_0^x
\,\sum_{n=0}^\infty  
         \frac{(-i\,t\,(x-\xi))^n}{n\,!} \, 
         g(\xi)\,d\xi\right) \, \widehat{f^+}(t) \, dt \\ 
         && \\
& = & \displaystyle \left(\frac{1}{2\,\pi} \, \int_{\mathbf R}
         {\mathcal J}^+ \, 
         \sum_{n=0}^\infty \left(-\,i\,t\,{\mathcal J}^+\right)^n \,
         \widehat{f^+}(t) \, dt \, g\right)(x) \\ 
         && \\
& = & \displaystyle \left(\left(\frac{1}{2\,\pi\,i} \int_{\mathbf R} {\mathcal
           J}^+/\left(1 + i\,t\,{\mathcal J}^+\right) \,
         \widehat{f^+}(t) 
         \, dt\right) \, g \, \right)(x) \\
&& \\
     & = & \displaystyle \lim_{\delta \to 0^+}\left(\frac{1}{2\,\pi\,i}\,\int_{\mathbf R}
      \widehat{f^+}(t)/\left(t\,I - i/\left({\mathcal
          J}^++\delta\,I\right)\right) \, dt \, 
      g\right)(x) \\ 
      && \\
      & = & \displaystyle \lim_{\delta \to 0^+}
      \left(\widehat{f^+}\left(i/\left({\mathcal J}^+ +
      \delta\,I\right)\right)\,g\right)(x) = 
      \left(\widehat{f^+}\left(i/{\mathcal J}^+\right)\,g\right)(x)\,,
  \end{array}
      \label{equation:2.12}
\end{equation} 

\noindent Similarly, $q^-(x) = \left(\widehat{f^-}\left(-i/{\mathcal
  J}^-\right)\,g\right)(x)$\,.  
  \end{Proof} $\blacksquare$

 \vspace{.2in} \noindent The above identities (2.7) are esoteric. They
 do, however readily enable applications. See \S4.6 of [Sp], or [Sc]. 

 \subsection{Fourier transform inversion.} 

We describe here an explicit novel formula for the inversion of 
Fourier integrals, namely for the determination of $f^\mp$
given $\widehat{f^\mp}$, where $\widehat{f^\mp}$ are the Fourier
transforms of $f^\mp$ as defined in (2.7).  
 
\newpage
\begin{Theorem}
  \label{theorem:2.6}
   Let ${\mathbf R}^+$ denote the interval $(0,\infty)$\,, and let us
   assume that we are given one or both of the functions
   $\widehat{f^\mp}$ where,

   \begin{equation}
        \widehat{f^\mp}(x) = \displaystyle \int_{{\mathbf R}^+} \,
        e^{\mp\,i\,x\,y} \, f^\mp(y)\,dy.
     \label{equation:2.13}
   \end{equation}

   \noindent Let the operators ${\mathcal J}^\mp$ have support on
   $(0,b) \subseteq {\mathbf R}^+$\,.  Then, for given
   $\widehat{f^\mp}$ on $(0,b)$\,,

    \begin{equation}
        f^\mp = \displaystyle \left(1/{\mathcal
          J}^\mp\right)\,\widehat{f^\mp}\left(\mp\,i/{\mathcal
          J}^\mp\right)\,1\,,
      \label{equation:2.14}
    \end{equation}

    \noindent where the ``$1$'' on the right hand side of (2.14) denotes the
    function that is identically $1$ at all points of $(0,b)$\,.

 \end{Theorem}

 \begin{Proof} We shall only prove this theorem for the case of
   $f^+$\,, inasmuch as the proof for the case of $f^-$ is similar.

 \vspace{.1in} \noindent Let us first recall, by assumption, that $f^\mp$ is
 differentiable.  Then
 
 \begin{equation}
   \begin{array}{rcl}
     f^+(x) - f^+(0) & = & \displaystyle \int_0^x \left(f^+(y)\right)^\prime\,dy
     \\
     && \\
     & = & \displaystyle \int_0^x \left(f^+\right)^\prime(x-y)\,dy \\
     && \\
     & = & \displaystyle \int_0^x \left(f^+\right)^\prime(x-y)\,1\,dy\,. 
   \end{array}
   \label{equation:2.15}
 \end{equation}

 \noindent This equation is now in the form of the equation for $q^+$
 in (2.9), except that $q^+$ in (2.9) is here replaced with the
 derivative $\left(f^+\right)^\prime$\,, and $g$ with the the function
 that has value $1$ on $(0,b)$\,.  Hence, using the Fourier
 transform $- i\,x\,\widehat{f^+}(x) + f(0)$ of
 $\left(f^+\right)^\prime$ and substituting into (2.15)\, we get the
 equation (2.10) for the case of $f^+$\,.
  \end{Proof} $\blacksquare$

 \vspace{.2in} \noindent The identities (2.10) and (2.15), while
 esoteric, can nevertheless yield applicable approximations in
 suitable analytic function settings, as they did for the case of
 Laplace transforms in [Sp] and in [Sc].

 \begin{Remark}
   \label{remark:2.7}

   If we substitute $y = \mp\,i\,\eta$ in (4.1), we are then back to
   the Laplace transform cases already covered extensively, starting
   with [Sp], \S4.5--4.6, and then followed up with all of the text,
   [Sc]\,. These sources cover the operator results of this section,
   elucidating them to approximation via use of Cauchy sequences of
   analytic functions based on Sinc methods of approximation.
   
   \vspace{.2in} \noindent We should add, that the two main theorems
   of this section have applications not only to approximation via
   Sinc methods, but also, to any other method of approximation,
   including methods that use orthogonal polynomials.

 \end{Remark}

 \section{Connection with interpolatory approximation}
 \setcounter{equation}{0}

 The formulas derived in the previous section are esoteric, but they
 have many applications when connected with interpolatory
 approximation\footnote{We could also include trigonometric
  polynomial in the examples which follow, e.g., those of [Sc], \S1.4,
  whose interpolate at points $x_j$ that are interior points of the
  interval of interpolation. Such formulas are effective for
  approximation of periodic functions.}.

\begin{description}
  \item{(i.)} For the case Legendre polynomials, the $x_j$ are the
  $n$ zeros of the Legendre polynomial $P_n(x)$ which are orthogonal
  on the interval $(a,b) = (-1,1)$ with respect to the weight function
  $w$ which is identically $1$ on $(-1,1)$\,;

  \item{(ii.)} For the case of Hermite polynomial interpolation, using
    the Hermite polynomials $H_n(x)$ that are orthogonal over
    ${\mathbf R}$ with respect to the weight function $w$\,, with
    $w(x) = \exp(-x^2)$, and with $H_n(x_j) = 0$ for $j =
    1\,,\ \ldots\,,\ n$\,; and

  \item{(iii.)} Other polynomials that are orthogonal with respect to a
    weight function, such as Jacobi polynomials, Gegenbauer
    polynomials, etc.
\end{description}

 \vspace{.2in} \noindent \begin{Definition}
\label{definition:3.1}
   Consider standard Lagrange interpolation at distinct points
 $x_j$\,, with $a < x_1 < x_2 \ldots x_n < b$

 \begin{equation}
   \begin{array}{rcl}
     f(x) & \approx & \displaystyle \sum_{k=1}^n \ell_k(x) \, f(x_k)\,, \\
     && \\
     \ell_k(x) & = & \displaystyle \prod_{j = 1\,,\ \ldots\,,\ n,\ \ j \neq
       k}\frac{x - x_j}{x_k - x_j}\,,
   \end{array}
   \label{equation:3.1}
\end{equation}

 \noindent we introduce a family of $n \times n$ matrices $A^\mp =
 \left[A_{j,k}^\mp\right]$ for which the
 entries $A_{j,k}^\mp$ are defined by

\begin{equation} 
\begin{array}{rcl}
A_{j,k}^+ & = & \displaystyle \int_a^{x_j} \ell_k(x) \, w(x) \, dx\,, \\
\\
&& \\
A_{j,k}^- & = & \displaystyle \int_{x_j}^b \ell_k(x) \, w(x) \, dx\,,
\end{array}
\label{equation:3.2}
\end{equation}

\noindent where $w$ is a weight function that is positive a.e. on
$(a,b)$\,, such that the moments $\int_a^b w(x)\,x^j\,dx$ exist for
every non--negative integer $j$\,.

\vspace{.1in} \noindent Setting

\begin{equation}
  \begin{array}{rcl}
    V\,f & = & (f(x_1)\,,\ \ldots\,,\ f(x_n))^T \\
    && \\
    L(x) & = & (\ell_1(x)\,,\ \ldots\,,\ \ell_n(x))\,,
  \end{array}
  \label{equation:3.3}
\end{equation}

\noindent and defining ${\mathcal J}^\mp_n$ by
  
\begin{equation}
  ({\mathcal J}^\mp_n\,w\,f)(x) = L(x) A^\mp\,V\,{\bf f}\,,
  \label{equation:3.4}
\end{equation}

\noindent so that if, for $D = \Im \mp z > 0$\,, and if $\hat{f}$ is
analytic in $D$\,, then the eigenvalues of $\mp i\,A^\mp$ lie in
$D$\,, so that the matrix $\hat{f}(\mp i\,(A^\mp)^{-1})$ is then well
defined, as also, is the approximation\footnote{Note that $V\,L =
  {\mathbf I}$\, with ${\mathbf I}$ the unit matrix, so that
  $V(\hat{f}\left(\mp\,i\,L\,(A^\mp)^{-1}\,V\right) =
  L(x)\,\hat{f}\left(\mp\,i (A^{\mp})^{-1}\right)\,V$\,.}

\begin{equation}
    V\,\hat{f}(\mp\,i/{\mathcal J}^\mp)\,g \approx \hat{f}\left(\mp
    i\,(A^\mp)^{-1}\right) \, V\,g\,.
  \label{equation:3.5} 
  \end{equation}

 \end{Definition}

 \vspace{.1in} \noindent
 \begin{Remark}
   \label{remark:3.2}
   \begin{description}
   \item{(i.)} The last line of (3.5) can be explicitly evaluated.  If
     for the case of $A^+$ (resp., if for the case of $A^-$ we have
     $A^+ = X\,\Lambda\,X^{-1}$ (resp., $A^- = Y\,\Lambda\,Y^{-1}$),
     where $\Lambda = {\rm diag}(\lambda_1,\ \ldots\,\ \lambda_n)$ is
     the diagonal matrix, which is the same for $A^+$ and $A^-$\,, and
     where $X$ (resp., $Y$) is the corresponding matrix of
     eigenvectors, then

 \begin{equation}
   \begin{array}{rcl}
     \widehat{f^+}\left(i\,(A^+)^{-1}\right) \, V\,g = X \,
     {\rm diag}(i/\lambda_1,\ \ldots,\ i/\lambda_n) \, X^{-1} \, V\,g\,,
   \end{array}
   \label{equation:3.6}
 \end{equation}

 \noindent and similarly for the term involving $\widehat{f^-}\,, A^-$
 and $Y$\,.   

 \item{(ii.)} If for (i.) above, the matrices $A^\mp$ are defined for
   $(-1,1)$\,, then for any other interval $(a,b)$\,, $A^\mp$ needs to
   be replaced with $C^\mp = (b-a)/2$\,; this means that the
   eigenvalues $\lambda_j$ are also to replaced with $(b-a)/2$\,, but
   the matrix of eigenvectors remains unchanged.
   \end{description}
   \end{Remark}
 
\vspace{.2in} \noindent {\bf Conjecture.}  We state the following
conjecture, for which the author of this paper offers \$300 for the
first proof or disproof:

\vspace{.1in} \noindent {\bf All of the eigenvalues of each the} $n
\times n$ {\bf matrices} $A^\mp$ {\bf that are defined as in
  Definition 3.1 for all polynomials} $\{p_n\}$ {\bf which are
  orthogonal over} $(a,b)$ {\bf with respect to the weight function}
$w$ {\bf lie on the open right half of the complex plane.}

\vspace{.1in} \noindent This conjecture has been shown to be true for
Sinc interpolation by Han and Xu [HX]; it has also been shown to be
true for Legendre polynomial interpolation by Gautschi and Hairer
[GH].  However, the conjecture as stated for Definition 3.1 is still
unproved for arbitrary weight functions $w$ that are positive a.e. on
$(a,b)$\,.

\vspace{.2in} \noindent The following result, can also be of use in
applications.  We select for this theorem the weighted Hilbert space
${\bf H}_w$ of all functions $f,\ g,\ \ldots,$ with inner product

\begin{equation}
  (f,g) = \displaystyle \int_{-1}^1 w(x) \, \overline{f(x)} \, g(x) \, dx\,,
  \label{equation:3.7}
\end{equation}

\noindent where $w$ is defined in Definition 3.1.  The
transformation of this formula to the interval $(a,b)$ via use of the
transformation $x = t(y) := (a+b)/2 + y \, (b-a)/2$ takes the form

\begin{equation}
  (f,g) = (F,G)_{(a,b)} =  \displaystyle \int_a^b
  W(y) \, \overline{F(y)} \, G(y) \, \frac{b-a}{2} \, dy\,,
  \label{equation:3.8}
\end{equation}

\noindent where $W(y) = w(t(y))$\,, $F(y) = f(t(y))$, and $G(y) =
g(t(y))$\,.  The proof of the following result is straight--forward,
and we omit it.

\newpage
\begin{Theorem}
  \label{theorem:3.3}
Let the operators ${\mathcal J}^\mp$ be defined as in Definition 2.2.
If $f \in {\mathbf H}_w$\,, then

\begin{equation}
  \begin{array}{l}
\displaystyle \left|\int_a^bW(y)\,\overline{F(y)} \, \int_a^y W(t) \,
F(t) \, dt\, dy \right| = \displaystyle \frac{(b-a)^2}{2} \,
\left|\int_{-1}^1 \overline{f(\eta)}\,\int_{-1}^\eta f(\xi) \, d\xi \,
d\eta\right| \,, \\
\\
\leq \displaystyle \frac{(b-a)^2}{2} \, \int_{-1}^1 \frac{d}{dx}
\int_{-1}^x w(t) \, |f(t)| \, dt \, w(x) \, |f(x)| \, dx \\
\\
\leq \displaystyle \frac{(b-a)^2}{2} \left(\int_{-1}^1
w(x)\,|f(x)|\,dx\right)^2\, dx\,.
\end{array}
  \label{equation:3.9}
\end{equation}

\end{Theorem}

\vspace{.2in} \noindent \section{Applications} 
 \setcounter{equation}{0}

We illustrate in this section, several examples that ensue by
approximation of the integration operators ${\mathcal J}^\mp$\,.  Such
approximations were first stated in [Sr], then proved by Kearfott [K],
Haber [H] and the author [Sp], \S4.5.  approximation based on using
the operators ${\mathcal J}^\mp$.  The use of these operators for
obtaining numerical solutions of differential and integral equations
was first discovered by the author in [Sp], and these were combined
extensively in [Sc], with formulas for approximating the indefinite
integral.  The present section illustrates applications by combining
indefinite integration with Lagrange polynomial approximation.

\subsection{Legendre polynomial approximation of a model}

\vspace{.2in} \noindent For sake of simplicity, we use only Legendre
polynomials as approximations, since applications using other bases
are dealt with in exactly the same way.  We also use polynomials of
degree at most $4$\,, in view of the present wide--spread interest in
problems arising from using large sets of data values.  Under the
assumption of analyticity of the reconstruction (see e.g., \S2 of
[Sc]), we get effective answers to problems whose solutions are
smooth, using a very small number of points to construct a solution.
That is, the explicit numerical examples of this section are most
effective for cases when the model to be approximated is smooth.

\vspace{.2in} \noindent It is well known that if a function $f$ is
analytic in a simply connected domain $D$ in the complex plane, and if
a closed interval $[a,b]$ is in the interior of $D$\,, then we can
approximate $f$ on $[a,b]$ via a polynomial of degree $n$ for which
the error approaches zero at a rate of ${\mathcal O}(\exp(- \, c\,
n))$ where $c$ is a positive constant.  To be more specific, if $(a,b)
= (-1,1)$, and if we let $\phi_j$ denote the normalized Legendre
polynomial of degree $j$\,, so that $\int_{-1}^1\phi_j(x) \, \phi_k(x)
\, dx = \delta_{j,k}$ where $\delta_{j,k}$ denotes the Kronecker
delta, then the numbers $c_k = \int_{-1}^1 f(x) \, \phi(x) dx$
approach zero exponentially, and the approximation

\begin{equation}
  f(x) \approx \sum_{k=0}^n c_k \, \phi_k(x)
\label{equation:4.1}
\end{equation}

\noindent will then be accurate for a relatively small value of $n$\,.  To
this end, we could start with a positive integer $m$ (with $m = 5$ in
this paper), then, letting $x_1\,,\ \ldots\,,\ x_m$ denote the
distinct zeros of $\phi_m$, and letting $w_1\,,\ \ldots\,,\ w_m$
denote the corresponding Gaussian integration weights, we could
evaluate numbers $c_1\,,\ \ldots\,,\ c_n$\,, with the formula

\begin{equation}
  c_k \approx \displaystyle \sum_{\ell=0}^m w_\ell\,f(x_\ell)\,\phi_j(x_\ell)
  \label{equation:4.2}
\end{equation}

\noindent and stop the process, when $c_k$ ($k \leq n$) is sufficiently
small, since $c_n$ is then of the order of the error of approximation.

\vspace{.2in} \noindent The above process can also be used when $f$ is
approximated via statistical sampling.  In this case, it would be
necessary to get a good approximation at the above defined zeros
$x_j$\,, and this can be done in many ways, one of which is
$\ell^1$--averaging.

\vspace{.1in} \noindent Note, also, the well known $m$--point Lagrange
interpolation formula, which is used extensively in this section, for
evaluating the polynomial $P_{m-1}$ which interpolates a function $f$
at the $m$ points $x_1\,,\ \ldots\,,\ x_m$:

\begin{equation}
  \begin{array}{rcl}
    P_{m-1}(x) & = & \displaystyle \sum_{j=1}^m \ell_k(x)\,f(x_k)\,,\\
    && \\
    \ell_k(x) & = & \displaystyle \prod_{j=1\,, j \neq k}^m \frac{x - x_j}{x_k
      - x_j}\,.
    \end{array}
  \label{equation:4.3}
\end{equation}

\vspace{.2in} \noindent
\subsection{Reconstruction from statistical data}

Statistical models can take on many different forms, see e.g., [X].
Typically, these can take the form of a system of nonlinear equations,
a system of ordinary differential equations (ODE) a system of partial
differential equations (PDE), or a system of integral equations, all
in the possible presence of noise. The variety of such equations means
that many different methods are in use for solving such equations.
Included among these are Fokker--Planck models [BS],  Navier--Stokes
equations [STB], Schr\"odinger equations [SKB] and the methods used
in [Sc] for solving partial differential and integral equations.
However, for sake of simplicity we shall restrict or presentation to
one--dimensional models, with approximation using Legendre
polynomials. 

\subsection{Exact formulas and their approximation}

\noindent The formulas which we shall describe below will initially be
defined in terms of the operators ${\mathcal J}^\mp$\,, and they thus
appear esoteric.  However, these operations can readily be
approximated via use of computable basis functions.  We shall only use
approximation via Legendre polynomial interpolation of degree 5 in
this section, although the programs can easily be altered to work for
arbitrary degrees, and with other bases -- for example, the text [Sc]
mainly uses sinc functions as bases -- although methods based on
Fourier polynomial bases are given in \S 1.4.7 and 3.10.2 of [Sc].
Explicit {\em Matlab} programs are available for all of the examples
of this section.

 \begin{description}
\item{(i.)} {\bf Fourier transform inversion.}  Statistical modeling at
  times -- such as in connection with vision [GW] -- involves the
  Fourier transforms.

 \vspace{.2in} \noindent We shall use the formula of Theorem 2 (a).

 \noindent Consider the trivial example for the recovery of $f$ on
 $(0,4)$ given $F^+(y) = \int_0^\infty f(t) e^{i\,y\,t} dt = 1/(1 -
 i\,y)$\,.  The exact solution is (*) $f(t) = e^{-t}$\,.  We use
 formulas (2.14 ) and (4.3), taking a Legendre polynomial of degree
 $5$, to get a matrix $A^+$ of order $5$ which we multiply by $2$ to
 get $C = 2\,A^+$, for approximation on $(0,4)$\,, twice the length of
 the interval $(-1,1)$\,.  By Theorem 2.6 and Equation (3.4), $y$ in
 the above Fourier transform is replaced with $i/C$\,, and also,
 selecting the new interpolation points $\xi_j = 2(1+x_j)$\,, we get

  \begin{equation}
   (f(\xi_1)\,,\ \ldots\,,\ f(\xi_n))^T \approx ({\bf I} + C)^{-1}\,{\bf
      1}\,,
    \label{equation:4.4}
 \end{equation}

 \noindent where ${\bf I}$ denotes the unit matrix, ${\bf 1}$ is a
 vector of $n$ ones, and where the matrix $({\bf I} + C)^{-1}$ is
 non--singular, since all of the eigenvalues of $A^+$ have positive
 real parts\footnote{Since $F$ is analytic in the upper half plane,
   and since the real parts of the eigenvalues of $A^\mp$ are
   positive, $F(i/A^+)$ is well defined for all such $A^+$\,, and so
   $({\bf I} + C)^{-1}$ is well defined.}.

 \vspace{.1in} \noindent Initially we just plot the five exact and
 computed values at $\xi_j = 1 + x_j$.  The next plots contain the
 exact solution $e^{-t}$\,, along with the computed approximation
 evaluated with the polynomial which interpolates the 5 computed
 values at $\xi_j$\,.  These plots are carried out on a $100$--point
 equi-spaced mesh.  Finally, we also plot the fine-mesh error, i.e.,
 the difference between the exact and computed solution at the same
 100 points.

 \vspace{0.2in}
 \begin{figure}[htb]
 \centering
\includegraphics[scale=0.75]{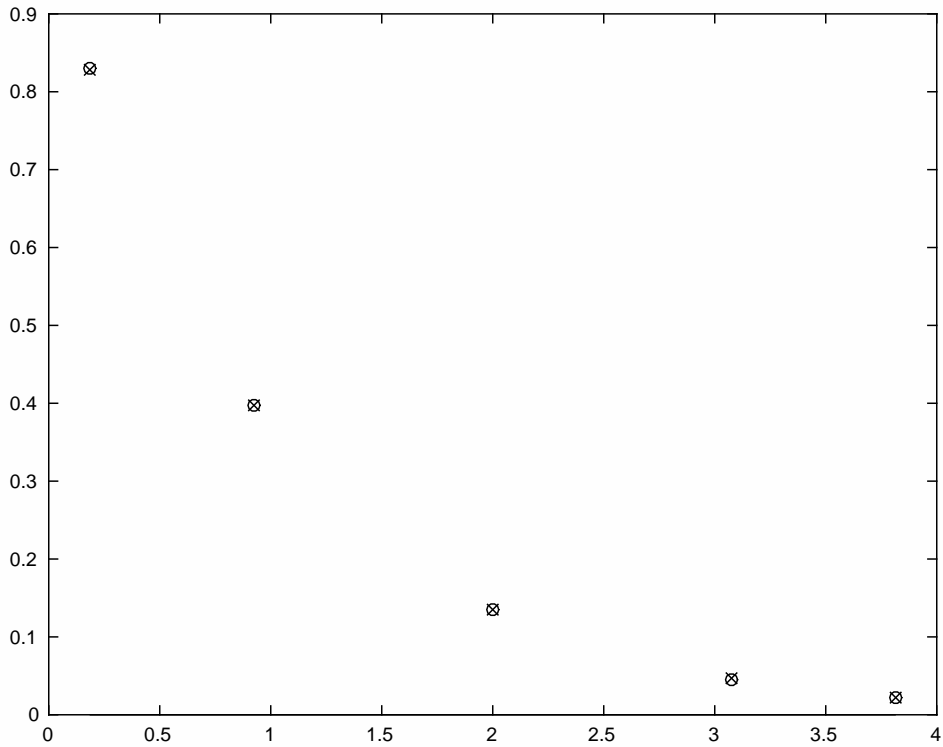}   
\caption{Course mesh plot of exact \& computed FT inversion}
\label{Figure:1}
 \end{figure}

 \vspace{0.2in}
 \begin{figure}[htb]
 \centering
\includegraphics[scale=0.75]{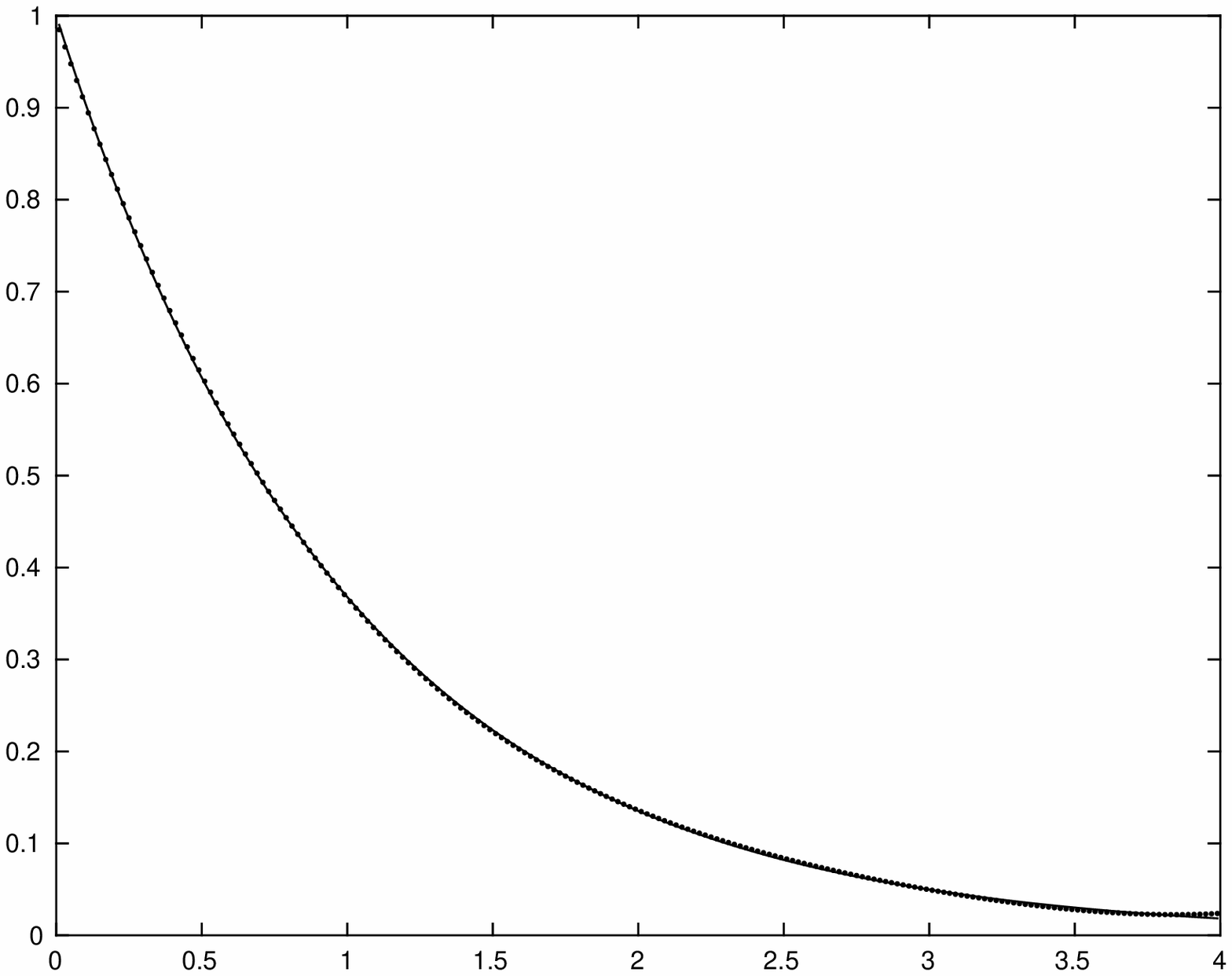}   
\caption{Fine mesh plot of exact \& computed FT inversion}
\label{Figure:2}
 \end{figure}
 
 \vspace{0.2in}
 \begin{figure}[htb]
 \centering
\includegraphics[scale=0.75]{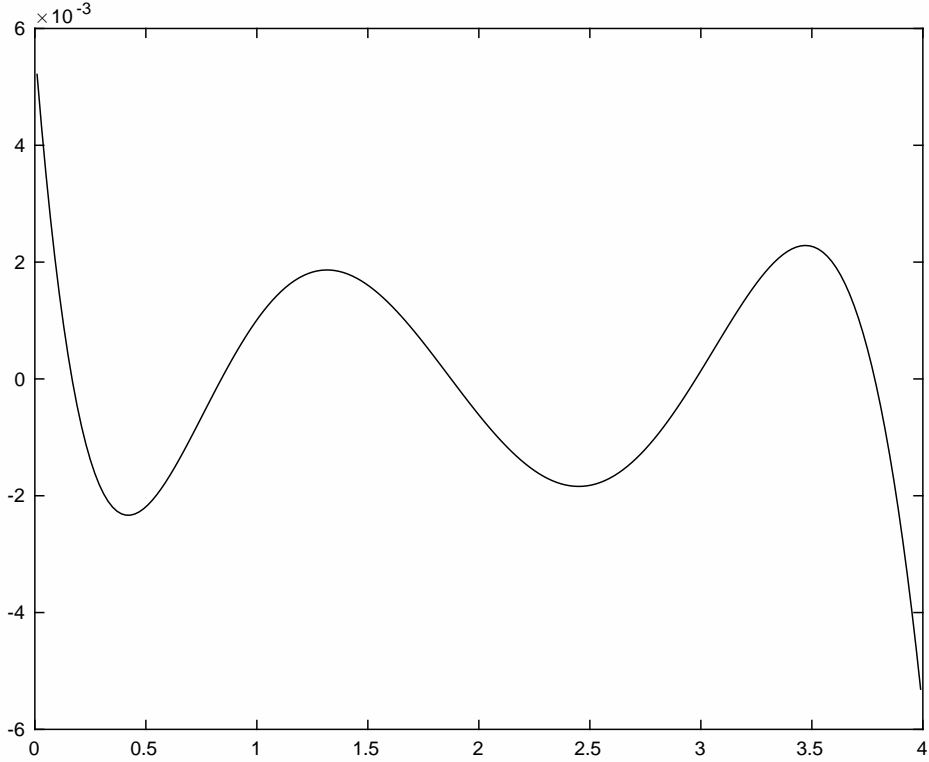}   
\caption{Fine mesh plot of error of FT inversion}
\label{Figure:3}
 \end{figure}

 \item{(ii.)} {\bf Laplace transform inversion}

  \vspace{.1in} \noindent The inversion formula for Laplace transform
  inversion was originally discovered by Stenger in [Sc].  The exact
  formula used here is only the third known exact formula for
  inverting the Laplace transform, the other two being due to Post [P]
  and Bromwich [B]., although we claim that the practical
  implementation of the Post formula is impossible, while the
  evaluation of the vertical line formula of Bromwich is both far more
  difficult and less accurate than our method, which follows.

  Consider the case of recovering the following function $f$ on $(0,2)$
  given $F$\,, where

\begin{equation}
\begin{array}{rcl}
f(t) & = & \displaystyle \frac{\sin(\pi\,t)}{\pi\,t};\\
&& \\
F(s) & = & \displaystyle \int_0^\infty f(t) e^{- st} dt = 1/2 -
(1/\pi) \tan^{-1}(s/\pi)\,.
\end{array}
\label{equation:4.5}
\end{equation}

\noindent We again use Lagrange polynomial approximation via
interpolation at the zeros of this polynomial of degree $5$\,.  The
length of the interval is $(0,2)$\,, which is the same as the length
of $(-1,1)$\,, for which the matrix $A^+$ is defined.  The new points
of interpolation are on the interval $(0,2)$ which shifts them by $1$
from the original interval $(-1,1)$ for Legendre polynomials, so that
the new points of interpolation are $\xi(j) = 1 + x(j)$\,.  The exact
inversion formula is $f = 1/{\mathcal J}^+ \, F(1/{\mathcal J}^+)\,1$\,,
with $F$ given on the right hand side of (4.5).  Hence replacing
${\mathcal J}^+$ with $A^+$\,, we get the following exact and
approximate solutions:

  \begin{equation}
    \begin{array}{l}
      f(\xi_j)  =  \displaystyle \frac{\sin(\pi\,\xi_j)}{\pi\,\xi_j}\,, \ \ \ j =
        1\,,\ \ldots\,,\ n\,, \\
        \\
        (f(\xi_1)\,,\ldots\,,\ f(\xi_n))^T \approx (A^+)^{-1}\,F((A^+)^{-1}) \,
        {\bf 1}\,,
    \end{array}
    \label{equation:4.6}
   \end{equation}
   
 \noindent where ${\bf 1}$ denotes a column vector of n
 ones\footnote{Note that the Laplace transform $F$ is analytic on the
   right half plane, and since the real parts of the eigenvalues of
   $A^+$ are positive, the matrix $F(A^+)$ is well defined.}.  In this
 case we need the eigenvalue representation $A^+ = X\,\Lambda\,X^{-1}$
 to evaluate last equation of (4.6).  We get

 \begin{equation}
   (f(\xi_1)\,,\ \ldots\,,\ f(\xi_n))^T \approx X\,\Lambda^{-1} \,
   F(\Lambda^{-1}) \, X^{-1} \, {\bf 1}\,.
   \label{equation:4.7}
 \end{equation}

 We next plot in fine mesh solution, which is obtained by evaluating
 the degree $4$ polynomial interpolation of the computed solution at
 100 equi--spaced points on $(0,4)$\,.  Finally, we also plot the
 fine-mesh error at the same 100 points.  The slight error at the
 end--points is due to the jump singularity of the Fourier transform at
 the origin.
 
\newpage
\begin{figure}[htb]
\centering
\includegraphics[scale=0.75]{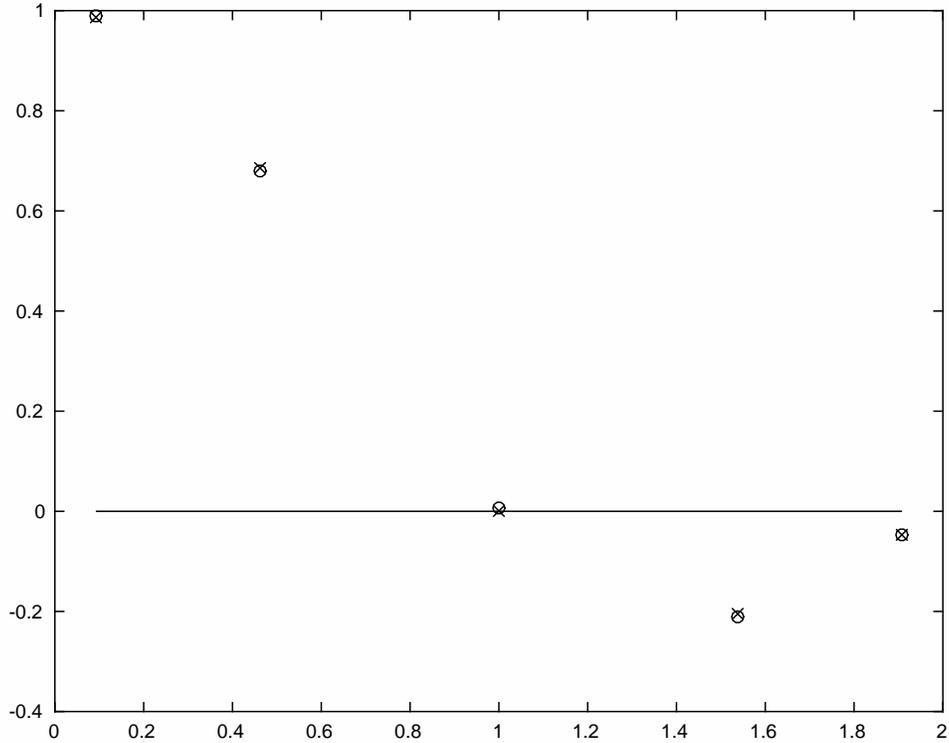}
\caption{Course mesh plot of exact (-) \& computed (.) LT inversion}
\label{Figure:4}
\end{figure}

 \vspace{0.2in}
\begin{figure}[htb]
\centering
\includegraphics[scale=0.75]{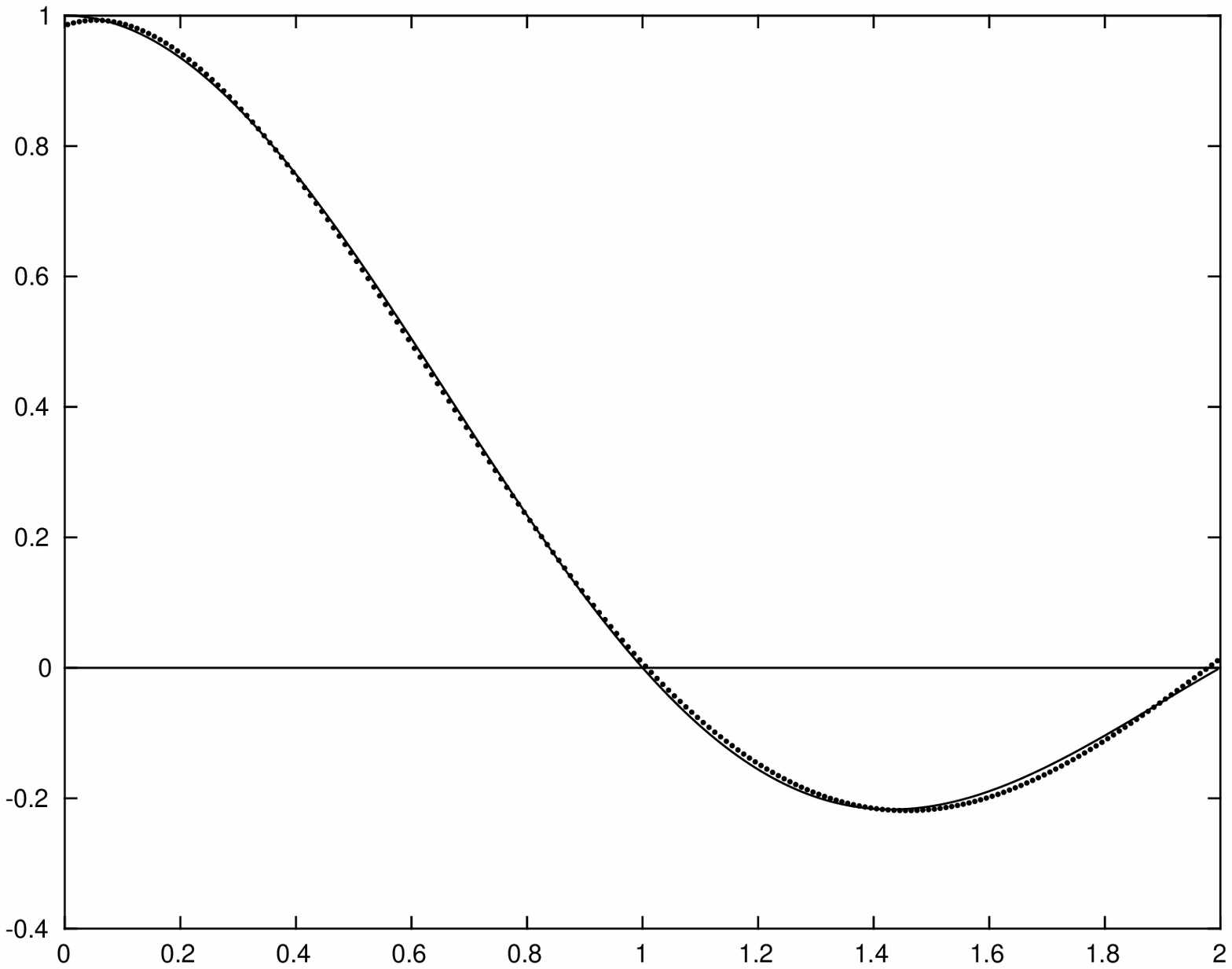}
\caption{Fine mesh plot of exact (-) \& computed (.) LT inversion}
\label{Figure:5}
\end{figure}

 \vspace{0.2in}
\begin{figure}[htb]
\centering
\includegraphics[scale=0.75]{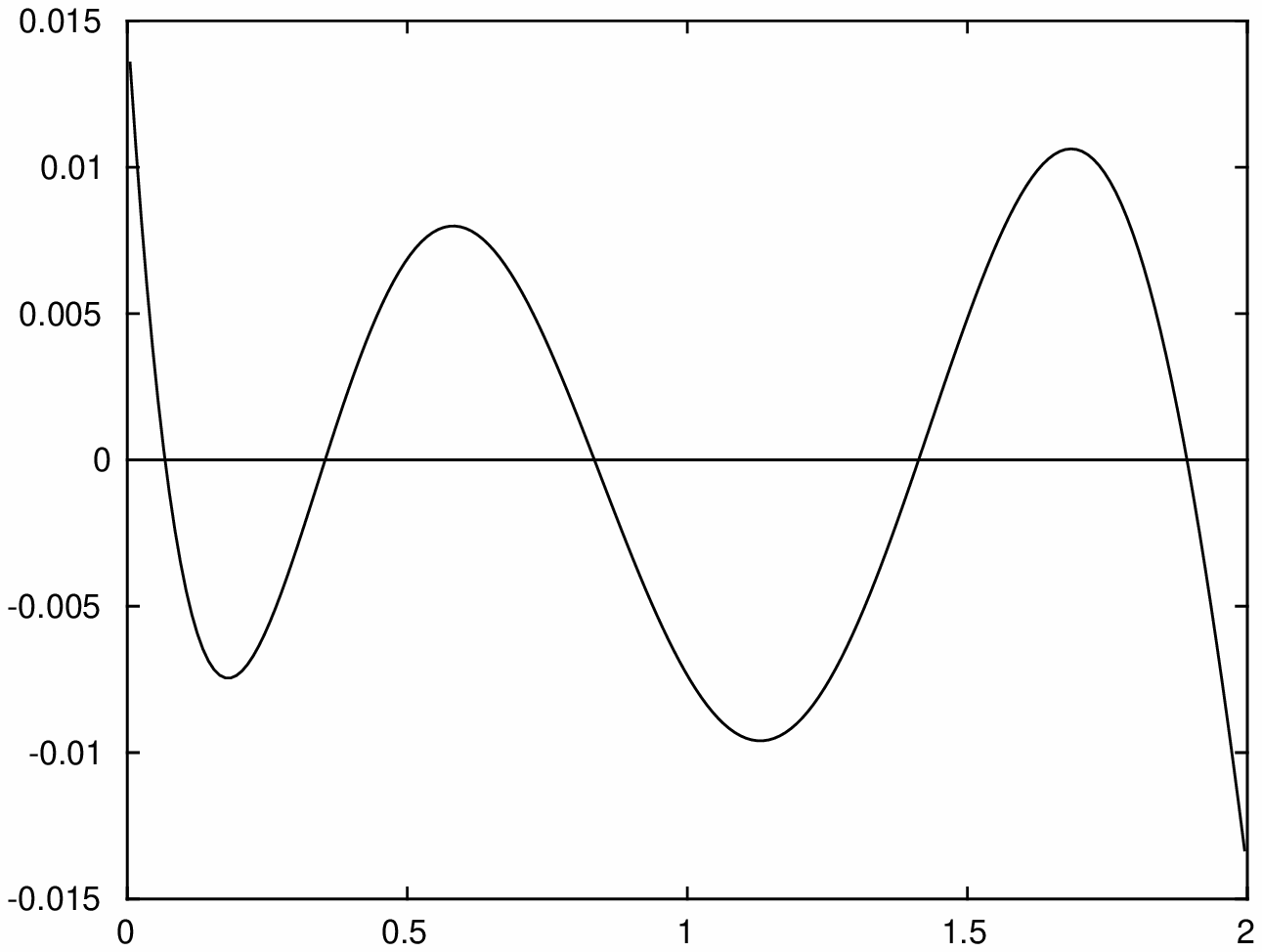}
\caption{Fine mesh plot of error of LT inversion}
\label{Figure:6}
\end{figure}

\item{(iii.)} {\bf Optimal control via Fourier transforms}

\noindent We illustrate here an application of Theorem 2.5.

\noindent Such an example may arise in simple design of a control, or
in the statistical determination of a feedback control, etc.  We wish
to evaluate integral

\begin{equation}
  p(t) = \int_t^3 f(t-\tau) \, g(\tau) \, d\tau\,,\ \ \ t \in (0,3)\,,  
  \label{equation:4.8}
\end{equation}

\noindent using the formula

\begin{equation}
   p = \widehat{f}(-i/{\mathcal J}^-) \, g\,, 
   \label{equation:4.9}
\end{equation}

\noindent where $\widehat{f}(y) = \int_{-\infty}^0 f(t) \, e^{i\,y\,t}
\, dt = \int_0^\infty f(-t) \, e^{-i\,y\,t} \, dt$\,.  Here we shall
use the matrix $A^-$ defined as in Definition 3.1 above, which must be
replaced by $C = (3/2)\,A^-$ for use on the interval $(0,3)$\,, and in
addition, the points of interpolation are $\xi_j = (3/2)\,(1 +
x_j)$\,.  Thus the approximation formula is

\begin{equation}
  V\,p \approx \widehat{f}(-i\,C^{-1}) \, V\,g\,.
  \label{equation:4.10}
\end{equation}

\vspace{.1in} \noindent We consider as an example, the evaluation of
the convolution integral

\begin{equation}
  \displaystyle \int_t^3 \exp(\alpha \, (t-\tau)) \, J_0(t-\tau) \,
  e^{-\,\beta\,\,\tau}\,d\tau\,, \ \ t \in (0,3)\,, 
  \label{equation:4.11}
\end{equation}

\noindent where $\alpha$ and $\beta$ are positive.  In this case, we
have the Fourier transform

\begin{equation}
  \widehat{f}(y) = \displaystyle \int_0^\infty e^{- i\,y\,t - a\,y\,t} \, J_0(t) \,
  dt =  \displaystyle \frac{1}{\alpha + i\,y} \, (1 + (\alpha + i\,y)^2)^{-1/2}\,,
  \label{equation:4.12}
\end{equation}

\noindent and upon replacing $y$ with $-i\,C^{-1}$\,, we get the
approximation

\begin{equation}
  \begin{array}{rcl}
      {\bf p} & := & (p(t_1)\,,\ \ldots\,,\ p(t_n))^T \\
   &&   \\
   & \approx &  C\,\left((1 + \alpha^2)\,C^2 +2 \, \alpha \,
      C \, + \, I\right)^{-1/2} \, {\bf g} \\
      && \\
      {\bf g} & = & \left(e^{-\,\beta\,t_1}\,,\ \ldots\,,\ e^{-\,\beta\,t_n}\right)^T\,.
  \end{array}
    \label{equation:4.13}
\end{equation}

\noindent The selection of several values of $\beta$ could be used to
model a given output $p$\,.  The example which follows uses $\beta =
0.7$\,.

\vspace{.1in} \noindent But more directly, since $p$ and $g$ are
related by the equation (4.12), we could also determine an accurate
control vector ${\mathbf g}$ to compute an accurate approximation to
the response:

\begin{equation}
{\bf g} \approx \left((1 + \alpha^2)\,C^2\, + 2\,\alpha\,C +
I\right)^{-1/2} \, {\bf p}\,.
\label{equation:4.14}
\end{equation}

\noindent Note here that the matrix multiplying $g$ in (4.13) can be
determined explicitly: setting $C = X\,\Lambda\,X^{-1}$\,, where
$\Lambda = {\rm diag}(\lambda_1\,,\ \ldots\,,\ \lambda_n)$\,, we have 

\begin{equation}
  \begin{array}{l}
D :=  X^{-1}\,\left((1 + \alpha^2)\,C^2\, + 2\,\alpha\,C
+I\right)^{1/2}\,X \,; \\
\\
D = {\rm diag}(d_1\,,\ \ldots\,,\ d_n); \ \ \ d_j = ((1 +
\alpha^2)\,\lambda_j^2 + 2\,\alpha\,\lambda_j + 1)^{1/2}\,.
\end{array}
\label{equation:4.15}
\end{equation}

\noindent This matrix is non--singular, and we can therefore use it to
compute an approximation to the function $g$ in (4.8), in order to 
achieve a any particular response.

\vspace{.2in} \noindent If we again use the same Legendre polynomial
of degree 5, we get a polynomial solution, but unfortunately, the
exact solution is not explicitly known.  To this end, we can use the
same equation (4.14), but with a larger value of $n$, to get a more
accurate solution; e.g., by taking $n = 11$, we get at least $8$
places of accuracy, which can be taken to be the exact answer for our
purposes, and which we use to compute the {\em exact} answer at $5$
points $\xi_j$\,.  We thus again get $3$ plots as above, i.e., as a
course mesh plot of the ``exact'' and our degree 4 approximation, a
fine mesh plot of the ``exact'' solution and of our degree 4
approximation of this solution, and a fine mesh plot of the difference
between these two quantities.

 \vspace{0.2in}
\begin{figure}[htb]
 \centering
\includegraphics[scale=0.75]{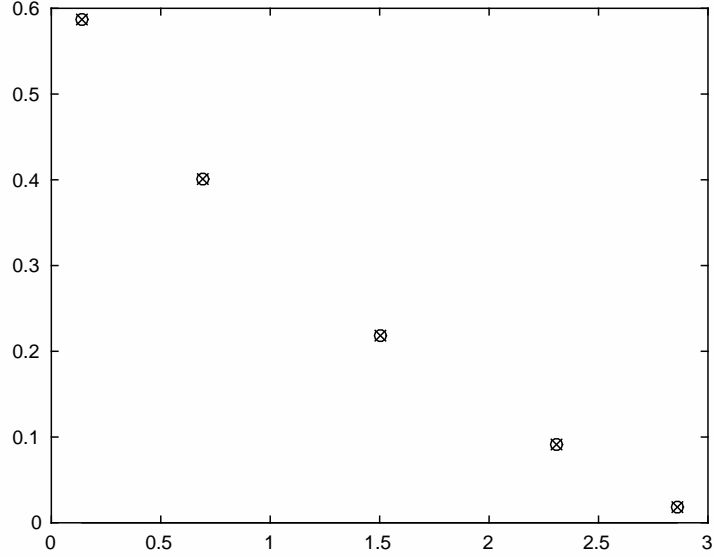}   
\caption{Course mesh plot of exact \& computed opt. control}
\label{Figure:7}
 \end{figure}

 \vspace{0.2in}
\begin{figure}[htb]
 \centering
\includegraphics[scale=0.75]{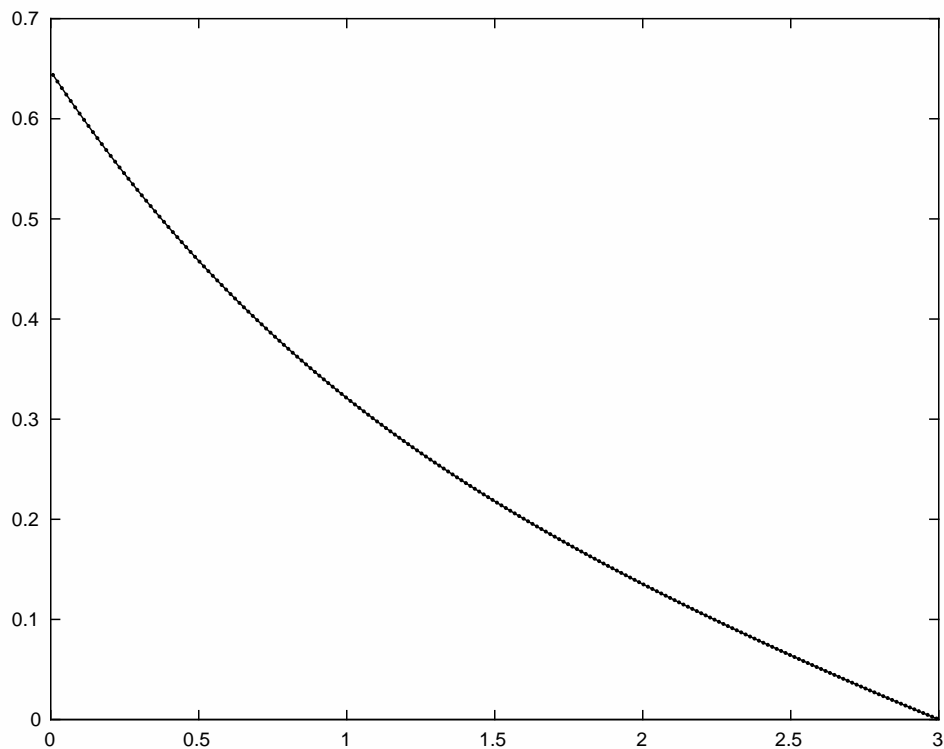}   
\caption{Fine mesh plot of exact \& computed opt. control}
\label{Figure:8}
 \end{figure}
 
 \vspace{0.2in}
\begin{figure}[htb]
 \centering
\includegraphics[scale=0.75]{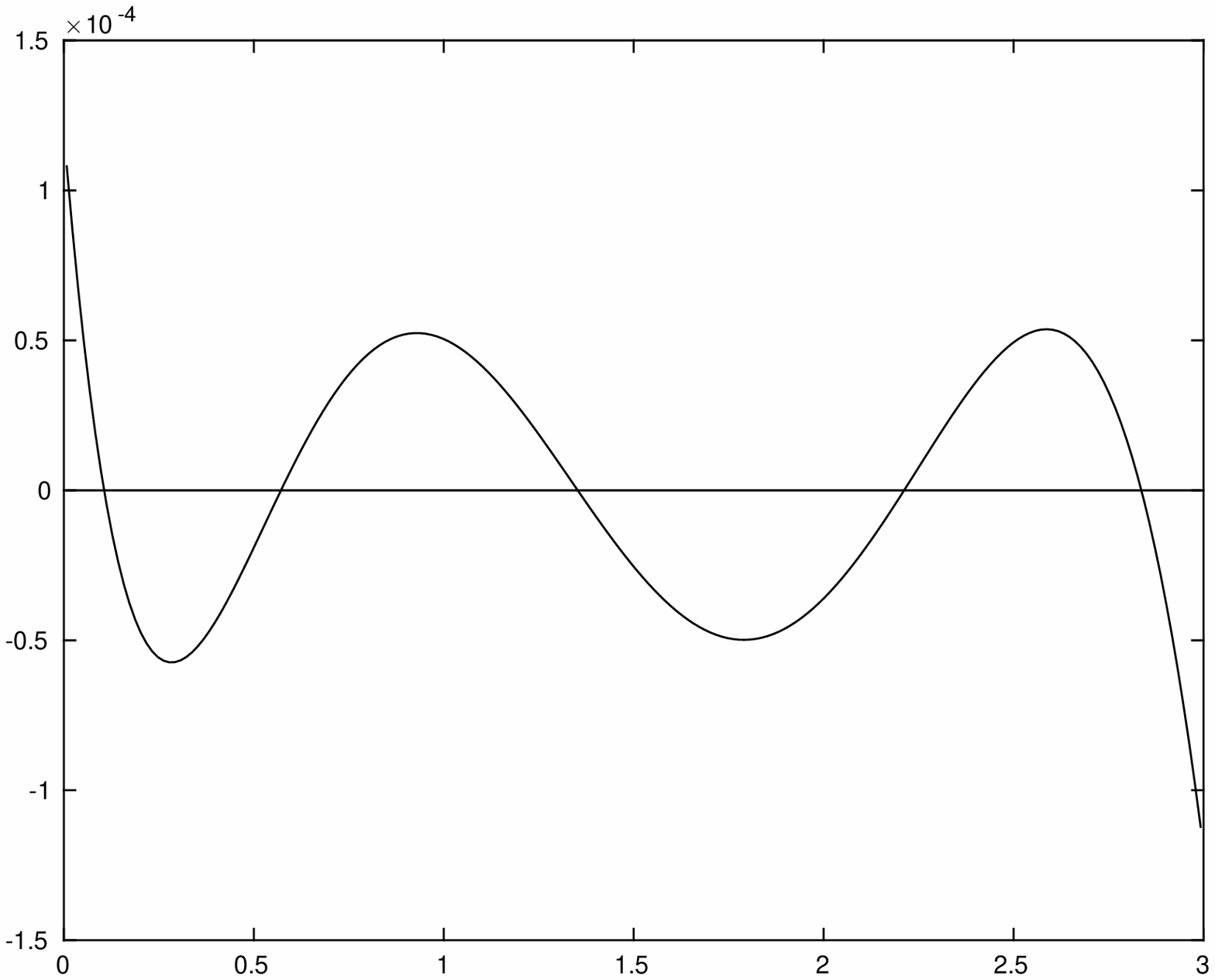}   
\caption{Fine mesh plot of error of opt. control}
\label{Figure:9}
 \end{figure}

\item {(iv.)} {\bf Modeling via ordinary differential equations}

  \noindent Most ODE solvers in use today are one step methods, and as
  such, their use is restricted because of stability, convergence,
  stiffness, and accuracy, and moreover, they are restricted to
  obtaining a solution on a finite interval.  Not so with the present
  method [DS], which extends its usage to polynomials the method of
  [SA] that was designed for Sinc approximation.

  \vspace{.1in} \noindent The most common ODE model for constructing
  methods of approximate solutions of ODE on an interval $(a,b)$ is

 \begin{equation}
   y^\prime = f(x,y)\,,\ \ \ y(a) = y_a \ \ \ ({\rm a\ constant})
\label{equation:4.16}
 \end{equation}

 \noindent Transforming to an equivalent integral equation, we get

 \begin{equation}
   y(x) = y_a + \displaystyle \int_a^x f(t,y(t))\,dt \ \ \ {\rm or}\ \ \ y =
   y_a + {\cal J}^+ \, f(\cdot,y(\cdot) \, 1\,.
   \label{equation:4.17}
 \end{equation}

 \vspace{.1in} \noindent Upon applying the approximation procedure of
 \S3 above, we can immediately convert The IE (4.17) to a system of
 algebraic equations

 \begin{equation}
   Y = y_a\,{\bf 1} + A^+\,{\bf f}\,{\bf 1}\,,
   \label{equation:4.18}
 \end{equation}

 \noindent where $Y = (y_1\,,\ \ldots\,,\ y_n)^T$\,, ${\bf I}$ is the
 unit matrix of size $n$\,, $A^+$ is as defined in \S3.2, and

 \begin{equation}
   {\bf f} := \left(\begin{array}{c}
     f(x_1,y_1) \\
     \vdots \\
     f(x_n,y_n)
   \end{array}
   \right)\,.
   \label{equation:4.19}
 \end{equation}

 \vspace{.1in} \noindent Some notes:
 \begin{description}
   \item{(1)}  Convergence will always occur for $b-a$ 
  sufficiently small, since it can be shown that the eigenvalues of 
  $\,A$ are bounded by $(b-a)/\sqrt{2}$\,.

\item{(2.)} Under a test of convergence criteria such as $\|Y^{(m)} -
  Y^{(m-1)}\| < {\varepsilon}$, the resulting solution will have
  polynomial degree $n$ accuracy at each of the points
  $z_1\,,\ \ldots\,,\ z_n$ where the interpolation is exact\,.  Under
  mild assumptions on $f$\,, we can then achieve similar accuracy at
  any set of points, e.g., at equi--spaced points.

\item{(3.)} Knowing accurate values of both $y$ and
  $y^\prime = f$ at $n$ points, we can get polynomial precision of
  degree $2\,n-1$  using Hermite interpolation.
\end{description}

\noindent Consider the case of recovering the following function $y$
on $(0,1/2)$\,, where

  \begin{equation}
    y\prime = 1 + y^2, \ \ \   y(0) = 0.
      \label{equation:4.20}
  \end{equation}

\noindent The exact solution is $y = \tan(t)$\,.

\vspace{.1in} \noindent Since the interval is $(0,1/2)$\,, the new
points of interpolation are $\xi_j = (1 + x_j)/2$\,, and the
corresponding indefinite integration matrix is $C = (1/2)\,A$\,.

By applying the above outlined method of approximation, we
replace the IE of (31) by the system of algebraic equations 

\begin{equation}
  Y = X + C\,Y_2\,,
  \label{equation:4.21}
\end{equation}

\noindent where $X = (t_1\,,\ \ldots\,,\ t_n)^T$\,, and $Y_2 =
(y_1^2\,,\ \ldots\,,\ y_n^2)^T$\,.  The interval $(0,1/2)$ was selected
here since, for example, if the matrix $C$ in (4.21) is replaced with
$A^+$, then Picard iteration does not converge to a solution of (4.20).
We could, of course, extend the solution further, by restarting it at
$t = 1/2$\,. 

\vspace{.1in} \noindent The computed solution and plots are the same
as in the previous cases, namely, a course plot of the exact and
computed solution at the points $\xi_j$\,, $j = 1\,,\ \ldots\,,\ 5$\,,
a fine mesh plot at $100$ equi--spaced points of the exact and the
computed polynomial solution on the interval $(0,1/2)$\,, and a fine
mesh plot of the difference between these functions on the same
interval.

 \vspace{0.2in}
\begin{figure}[htb]
 \centering
\includegraphics[scale=0.75]{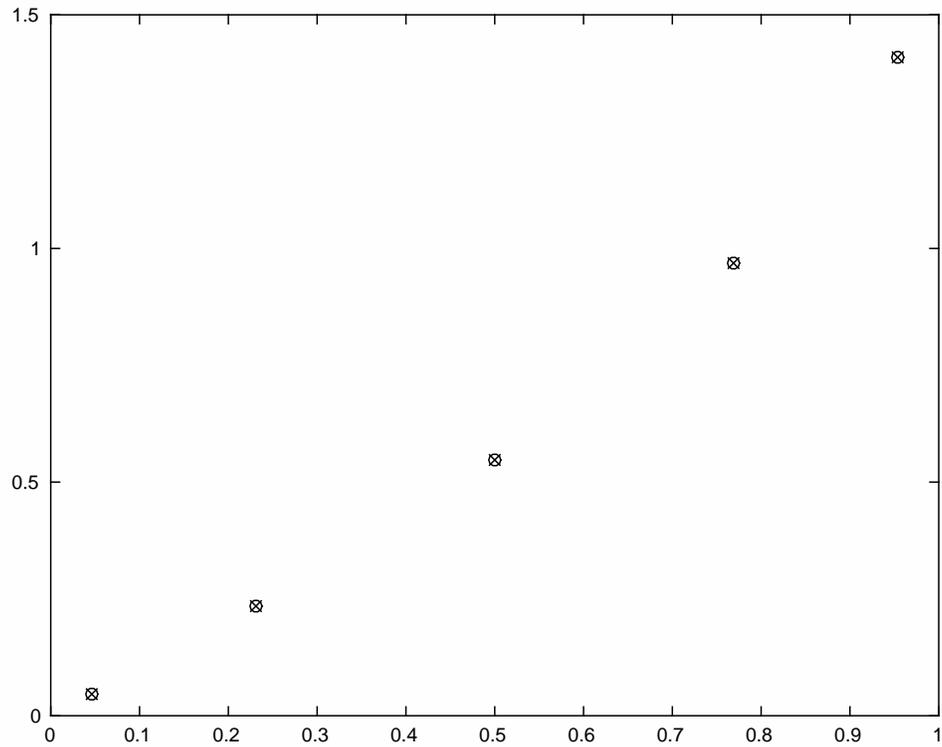}   
\caption{Course mesh plot of exact \& computed ODE}
\label{Figure:7}
 \end{figure}

 \vspace{0.2in}
\begin{figure}[htb]
 \centering
\includegraphics[scale=0.75]{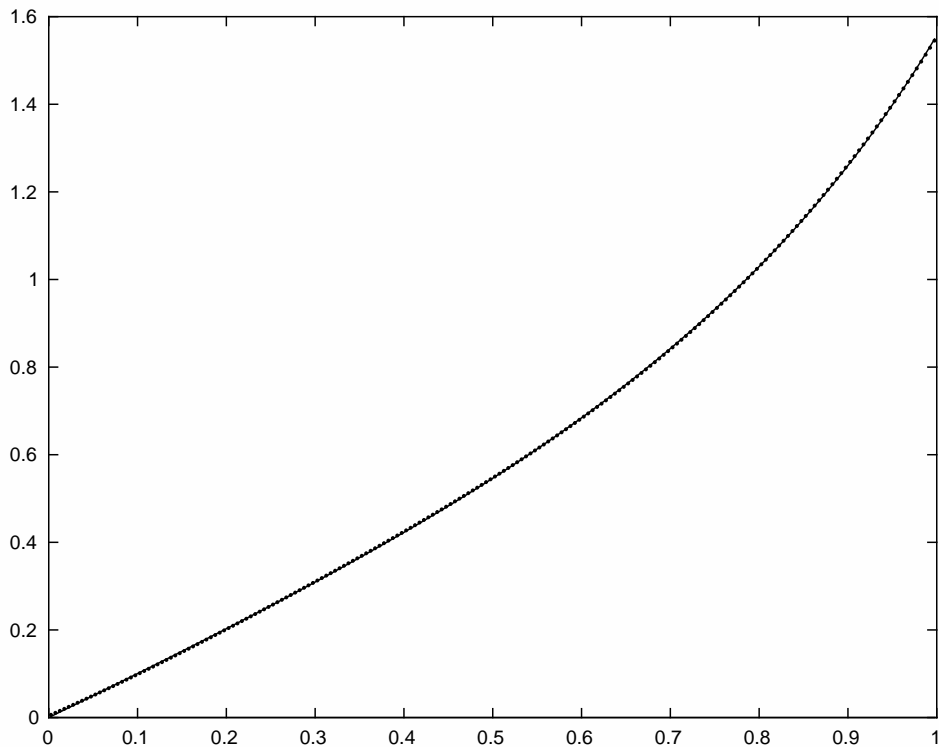}   
\caption{Fine mesh plot of exact \& computed ODE}
\label{Figure:8}
 \end{figure}
 
 \vspace{0.2in}
\begin{figure}[htb]
 \centering
\includegraphics[scale=0.75]{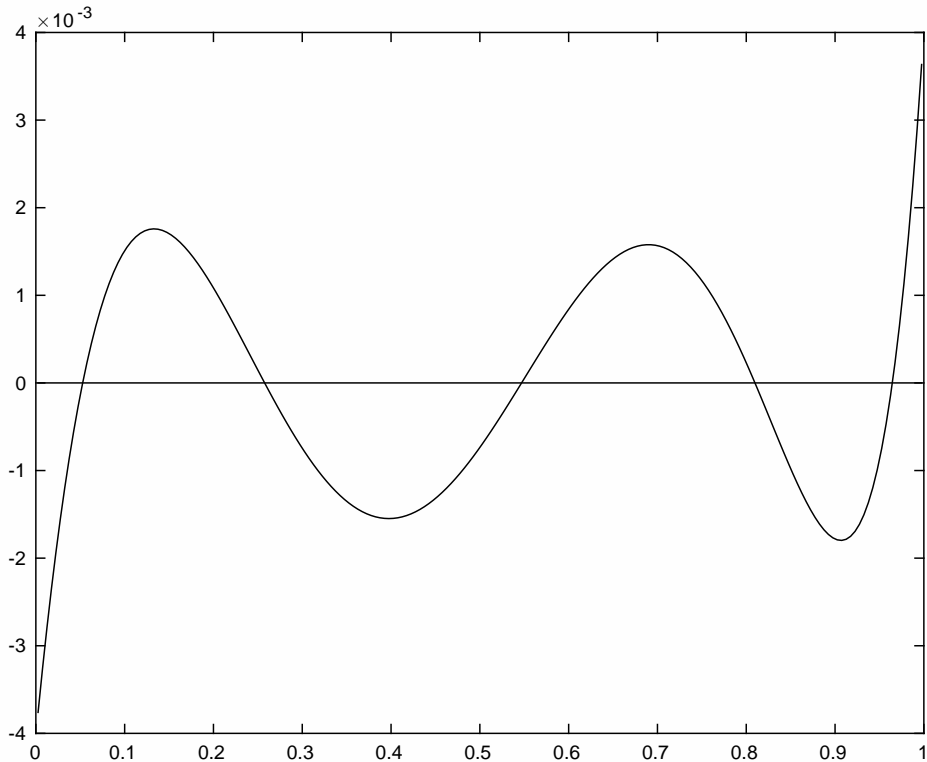}   
\caption{Fine mesh plot of error of ODE}
\label{Figure:9}
 \end{figure}

\item{(v.) Modeling via Wiener--Hopf Equations}

\noindent The classical Wiener--Hopf integral equation with solution
$f$ for given $k$ defined on ${\bf R}$ and $g$ defined on $(0,\infty)$
takes the form

\begin{equation}
  f(x) - \displaystyle \int_0^\infty k(x - t) \, f(t) \, dt = g(x)\,,\ \ x \in
  (0,\infty)\,.
  \label{equation:4.22}
\end{equation}

Many thousands of papers have been written on the solution of this
equation, particularly with reference to the mathematically beautiful
factorization procedure originally discovered by Wiener and Hopf in
1931 for solving this equation.  Unfortunately, such a factorization
cannot be determined for nearly all problems\footnote{An explicit
  factorization is, in fact, known, for the case when $k \in {\bf
    L}^1({\bf R})$ and $g \in {\bf L}^1(0,\infty)$\,, but this does
  not lend itself to a practically efficient method.} of the type (4.22).
A revision to this mathematics among the top 10 mathematics departments
occurred from about 1960 until 1970, but none of this pure mathematics
activity provided any insight to the solution of (4.22).  We illustrate
here a method of solving this problem, and while our illustration is
almost trivial, in that the exact equation is easier to solve than our
approximating equation, our approximating method nevertheless
applies to all equations of the the type (4.22).  Moreover, the solution
we present in what follows is both efficient and accurate.

\vspace{.1in} \noindent By splitting the definite convolution integral
in (4.22) as an integral from $0$ to $x$ plus an integral from $x$ to
$\infty$\,, to get two indefinite integrals, so that (4.22) can be
rewritten as

\begin{equation}
  f(x) - \displaystyle \int_0^x k(x-t)\,f(t)\,\,dt - \int_x^\infty
  k(x-t)\,f(t)\,dt = g(x)\,, \ \  x > 0\,.
  \label{equation:4.23}
\end{equation}

\noindent At this point we can invoke Theorem 2.5, which enables us to
replace the convolution integrals.  Letting $\widehat{k^+}$ denote the
Fourier transforms of $k$ taken $(0,\infty)$ (resp., letting
$\widehat{k^-}$ denote the Fourier transform of $k$ taken over
$(-\infty,0)$) we get the ``exact solution'',  

\begin{equation}
  f - \left(\widehat{k^+}\left(i \, (J^+)^{-1}\right)  - 
  \widehat{k^-}\left(-i \, (J^-)^{-1}\right) \, \right) f  = g\,,
   \label{equation:4.24}
 \end{equation}

 \noindent or, in collocated form, and now using both matrices
 $A^\mp$\,, we get 

 \begin{equation}
   \begin{array}{l}
     (f_1\,,\ \ldots\,,\ f_n)^T \approx \\
 \ \ \     \displaystyle \left({\bf I} - 
   \widehat{k^+}\left(i\,(A^+)^{-1}\right) - 
   \widehat{k^-}\left(-i\,(A^-)^{-1}\right)\right)^{-1} \,
   (g_1\,,\ \ldots\,,\ g_n)^T\,.
   \end{array}
 \label{equation:4.25}
 \end{equation}

\noindent In (4.25), ${\bf I}$ denotes the unit matrix, while the subscripts
``${(\cdot)_j}$'' on $f$ and $g$ denote values to be computed and exact
values respectively.   It should be observed that this representation
also yields accurate approximate solutions in cases when (4.22) has
non-unique solutions; such sulutions can be obtained by using singular
value decomposition to solve (4.25). 

\vspace{.2in} \noindent For example, consider the (rather trivial to
solve) equation,

 \begin{equation}
       f(x) - \displaystyle \int_0^1 k(x - t) f(t) dt = g(x) , \ \ \ 0
       < x < 1 .
  \label{equation:4.26}      
\end{equation}

\noindent where

  \begin{equation}
 \begin{array}{rcl}
 k(x) = \left\{ \begin{array}{rcl}
 - e^{-\,x} & {\rm if} & -1 < x < 0\,, \\
 - e^{-x} & {\rm if} & x > 0\,. \\
 0 & {\rm if} & -\infty < x < -1.
 \end{array}
 \right.
 \end{array}
  \label{equation:4.27}
  \end{equation}

 \noindent This equation is easier to solve analytically than
 numerically, although this is not the case for most functions $k$\,.
 We could just as easily solve (4.22) for a more complicated kernel
 such as $k(x) = \log(x)\,e^{- x}/(1+x^2)^{0.3}$ on $(0\,,\infty)$\,,
 and with $k(x)$ having a different, but a similarly complicated
 expression for $x < 0$\,.  But the procedure is the same in both the
 more complicated case as for this case.  The equation has the unique
 solution, $f(t) = g(t) - c\,e^{-t}$\,, with $c =
 (1/2)\,\int_0^1\,e^t\,g(t)\,dt$\,, for any $g$ defined on
 $(0\,,1)$\,.  In particular, if $g(t) = 2\,e^{-1/2}\,t\,e^{t^2 - t}$
 then $c = \sinh(1/2)$\,. In this case, we have\footnote{Notice, the
   integral for $\widehat{k}_-$ must be truncated, since it will not
   converge otherwise.  In particular, the integration from $-2$ to
   $0$ instead of from $-1$ to $0$ served to avoid the singularity due
   to the truncation of $\widehat{f}_-$ at $-1$\,.}

  \begin{equation} 
 \begin{array}{rcl}
 \widehat{k}^+(x) & = & - \displaystyle \int_{{\mathbf R}^+}  e^{i\,x\,y} \,
 e^{- y} \, dy = - 1/(1-i\,x) \\
 && \\
\widehat{k}^-(x) & = & - \displaystyle \int_{-2}^0 e^{i\,x\,t} \,
e^{-t} \, dt \\
& = & \displaystyle 2\,\exp(1 - i\,x) \, \frac{\sinh(1 - i\,x)}{1 - i\,x}\,.
 \end{array}
  \label{equation:4.28}
 \end{equation}

  \noindent Hence the operators $\widehat{k}_+(i\,({\cal J}^+)^{-1})$
  and $\widehat{k}_-(-i\,({\mathcal J}^-)^{-1})$ can be explicitly
  expressed, and replacing ${\mathcal J}^+$ with $A^+$ and ${\mathcal
    J}^-$ with $A^-$\,, we get

 \begin{equation} 
 \begin{array}{l}
 \widehat{k}^+(i\,(A^+)^{-1}) = - A^+\,(I + A^+)^{-1}\,, \\
  \\
  \widehat{k}^-(-i\,(A^-)^{-1}) = 2\,\exp\left(I -
  (-i\,(A^-)^{-1})\right)  \,\cdot \, \\
  \ \ \cdot \, \sinh\left(I
    - i(-i(A^-)^{-1})\right)\,\left(I - i\,-i(A^-)\right)^{-1} .
     \end{array}
 \label{equation:4.29}
  \end{equation}

 \noindent These matrices can be readily 
 computed, e.g., if $\lambda_j$ ($j = 1,\ \ldots,\ n$) are the
 eigenvalues of $A^\mp$\, then, setting $u_j = -
 (1+1/\lambda_j)$\,,\ \ \  $v_j = - (1 - 1/\lambda_j)$\,, and $w_j =
 \exp(-v_j)\,\sinh(v_j)/v_j$\,, $U = {\rm diag}(u_1,\ \ldots,\ u_n)$
 and $W = {\rm diag}(w_1,\ \ldots,\ w_n)$\,, we get
 $\widehat{k^-}\left(-i\,(A^-)^{-1}\right) = Y\,U\,Y^{-1}$\,, and
 $\widehat{k^-}\left(-i\,(A^+)^{-1}\right) = X\,W\,X^{-1}$\,.   

 \noindent Figure 1 is a plot of the exact and the computed solution,
 where to get the approximate solution, by using a $4$--degree
 polynomial approximation on $(0,1)$\,.  We could, of course easily
 have gotten greater accuracy using a higher degree Legendre
 polynomial.  As above, we have again plotted the course mesh exact
 and approximate solution, the fine mesh approximate and the fine mesh
 of the difference between the exact and the computed solution.

 \vspace{0.2in}
 \begin{figure}[htb]
 \centering
\includegraphics[scale=0.75]{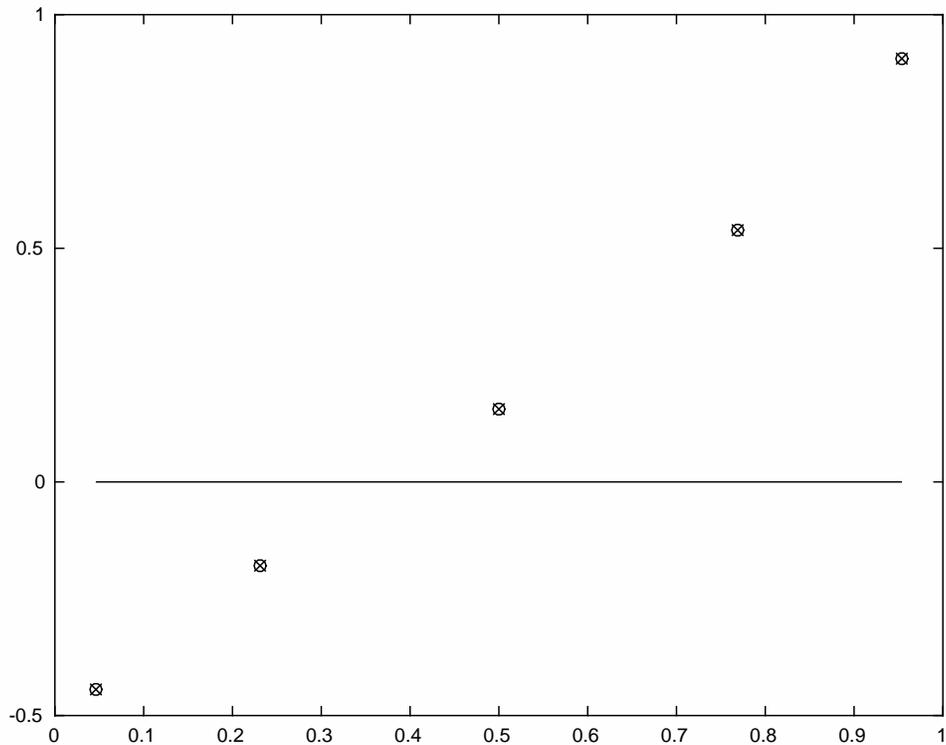}   
\caption{Course mesh plot of exact \& computed Wiener--Hopf}
\label{Figure:10}
 \end{figure}

 \vspace{0.2in}
 \begin{figure}[htb]
 \centering
\includegraphics[scale=0.75]{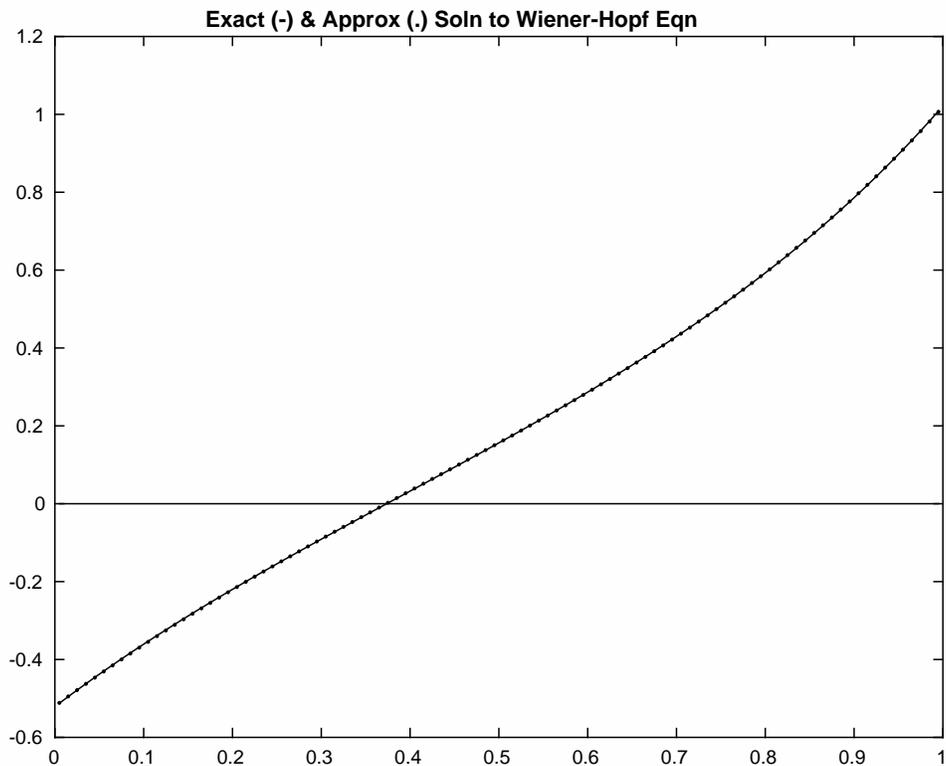}   
\caption{Fine mesh plot of exact \& computed Wiener--Hopf}
\label{Figure:11}
 \end{figure}
 
 \vspace{0.2in}
 \begin{figure}[htb]
 \centering
\includegraphics[scale=0.75]{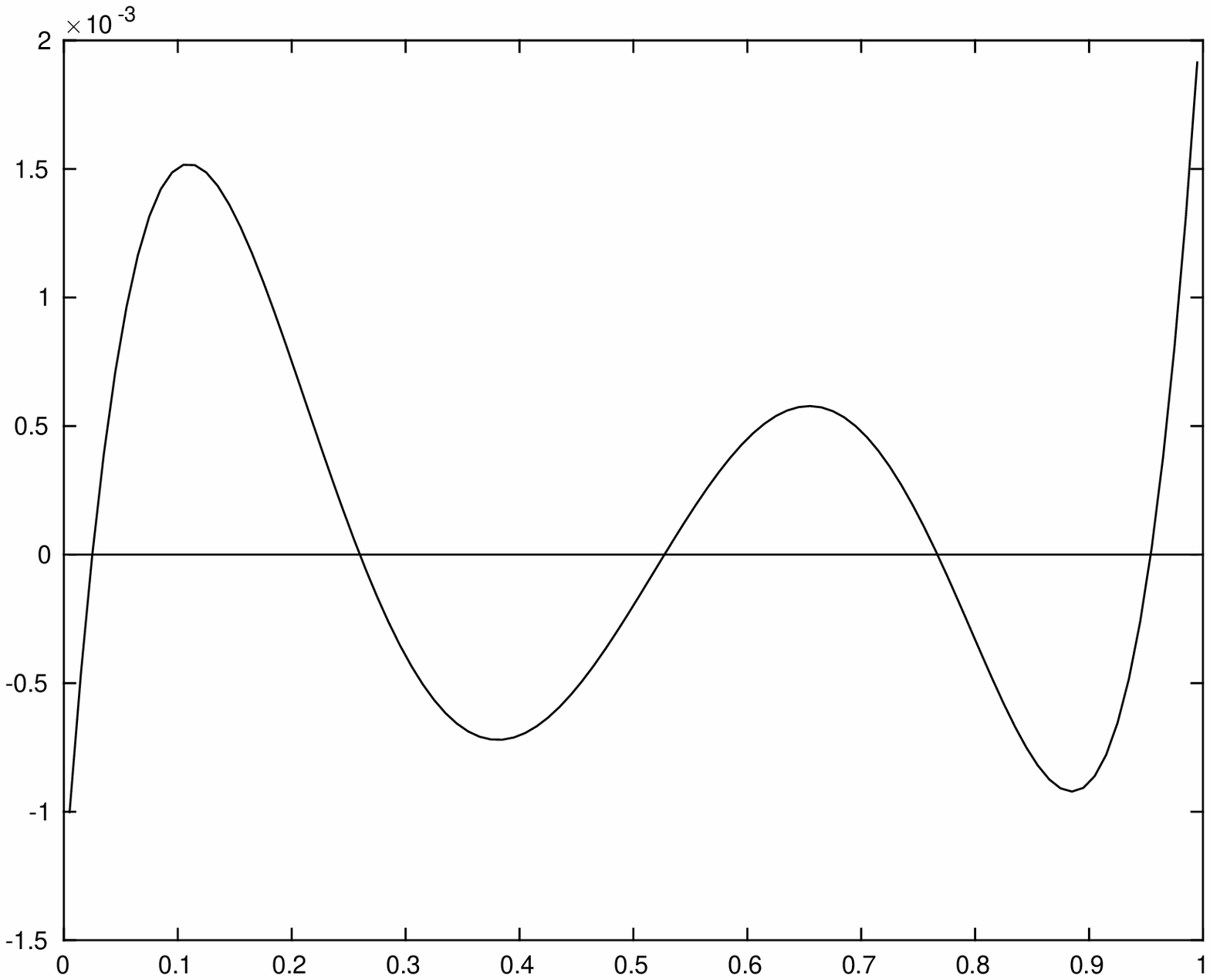}   
\caption{Fine mesh plot of error of Wiener--Hopf}
\label{Figure:12}
 \end{figure}
   
  \end{description}

\section*{References}

\begin{description}
  \begin{sloppypar}
\item{[B]}  {\bf T. Bromwich}, {\em Normal coordinates in dynamical 
  systems}, Proc. London Math. Soc., {\bf 15} (1916) 401--448. 
\item{[BS]} {\bf G. Baumann \& F. Stenger}, {\em Fractional Fokker--Planck
  Equation}, Mathematics {\bf 5} (2017) 1--19.
\item{[DS]}  {\bf G. Dahlquist} \& {\bf F. Stenger},  {\em
  Approximate solution of ODE via approximate indefinite integration},
  submitted.
\item{[GH]} {\bf W. Gautschi \& E. Hairer}, {\em A conjecture of Stenger in
  the theory of orthogonal polynomials}, submitted.
\item{[GR]} {\bf K.E. Gustafson \& D.K.M. Rao}, {\em Numerical Range:
     The Field of Values of Linear Operators and Matrices},
     Springer--Verlag (1996).  
\item {[H]} {\bf S. Haber}, {\em Two Formulas for Numerical Indefinite
  Integration}, Math. Comp., v.  60 (1993) 279­296.
\item{[HX]}  {\bf L. Han \& J. Xu}, {\em Proof of Stenger's Conjecture on Matrix
  $I^{(-1)}$ of Sinc Methods}, J. Comp. Appl. Math., {\bf 255} (2014)
  805--811.
 \item{[K]} {\bf R.B. Kearfott}, {\em A Sinc Approximation for the
   Indefinite Integral}, Math. Comp., V. 41 (1983) 559--572.
\item{[L]}  {\bf P. Lax}, {\em Functional Analysis}, Wiley \& Sons (2002).
\item{[Fc]} {\bf F. Stenger}, {\em Collocating Convolutions,}  Math. Comp.,
  {\bf 64} (1995) 211--235.
\item{[N]} {\bf A. Naghsh--Nilchi}, {\em Sinc Convolution Method of Computing
      Solutions to Maxwell's Equations}, Ph.D. thesis (1997).
  quarter, 1997.
\item{[P]} {\bf E. Post}, {\em Generalized differentiation},
  Trans. AMS {\bf 32} (1930) 723-781.
\item {[SA]} {\bf F. Stenger, S\AA. Gustafson, B. Keyes,
  M. O'Reilly, \& K.  Parker}, {\em ODE -- IVP -- PACK via Sinc
  Indefinite Integration and Newton's Method}, Numerical Algorithms
  {\bf 20} (1999) 241--268.
\item{[Si]} {\bf F. Stenger}, {\em Numerical Methods Based on the Whittaker
  Cardinal, or Sinc Functions}, SIAM Rev. 23 (1981) 165--224.
\item{[Sp]} {\bf F. Stenger}, {\em Numerical methods based on Sinc and
   analytic functions}, Springer--Verlag (1993).
\item{[Sc]} {\bf F. Stenger}, {\em Handbook of Sinc numerical methods},
  CRC Press (2011).
\item{[Srh]} {\bf F. Stenger}, {\em A proof of the Riemann
  hypothesis}, \\ http://arxiv.org/abs/1708.01209
\item{[STB]} {\bf F. Stenger, D. Tucker \& G. Baumann},  {\em Solution of
  Navier--Stokes on ${\mathbf R}^3 \times (0,T)$}, Springer--Verlag (2016).
\item{[SKB]} {\bf F. Stenger, V. Koures \& K. Baumann}, {\em Computational
    Methods for Chemistry and Physics, and Schr\"odinger $3 + 1$}, In:
  John R. Sabin and Remigio Cabrera-Trujillo, editors, Advances in
  Quantum Chemistry, Burlington: Academic Press, {\bf 71} (2015)
  265-298.
\item{[Y]} {\bf T. Yamamoto}, {\em Approximation of the Hilbert transform
   via use of Sinc convolution}, ETNA  {\bf 23} (2006) 320--328. 
  \end{sloppypar}
\end{description}

\end{document}